\documentstyle[12pt]{article}

\newtheorem{theorem}{Theorem}[section]     \newtheorem{lemma}[theorem]{Lemma}
\newtheorem{corollary}[theorem]{Corollary} \newtheorem{claim}[theorem]{Claim}
\newtheorem{problem}{Problem}
\newcommand{\subsec}[1]{\par\bigskip\begin{center}{{\bf #1}}\end{center}}
\newcommand{\pf}{\par\noindent{\bf Proof:}\par}
\newcommand{\qed}{\nopagebreak\par\noindent\nopagebreak$\blacksquare$\par}
\newcommand{\sone}{{\sf S}_1}     \newcommand{\sfin}{{\sf S}_{fin}}
\newcommand{\zz}{{\Bbb Z}^\omega} \newcommand{\ufin}{{\sf U}_{fin}}
\newcommand{\op}{{\cal O}}        \newcommand{\om}{\Omega}
\newcommand{\lm}{\Lambda}         \newcommand{\ga}{\Gamma}
\newcommand{\la}{\langle}         \newcommand{\ra}{\rangle}
\newcommand{\oo}{\omega^\omega}   \newcommand{\non}{{\sf non}}
\newcommand{\split}{{\sf Split}}  \newcommand{\covmeag}{{\sf cov}({\cal M})}
\input amssym.def
\input amssym

\begin{document}
\begin{center}
{\large \bf The combinatorics of open covers (II)}
\end{center}

\bigskip

\begin{center}
Winfried Just\footnote{partially supported by NSF grant
DMS-9312363},
Arnold W. Miller,
Marion Scheepers\footnote{partially supported by NSF grant DMS-95-05375},\\
and Paul J. Szeptycki
\end{center}

\bigskip

\begin{abstract} We continue to investigate various
 diagonalization properties
 for sequences of open covers of separable metrizable spaces
 introduced in Part I.  These properties generalize
 classical ones of Rothberger, Menger, Hurewicz, and
 Gerlits-Nagy.  In particular, we show that most of the properties
 introduced in Part I are indeed distinct. We characterize
 two of the new properties by showing that they are equivalent
 to saying all finite powers have one of the classical properties
 above (Rothberger property in one case and in the Menger property
 in other).  We consider for each property the smallest
 cardinality of metric space which fails to have that property.
 In each case this cardinal turns out to equal another well-known
 cardinal less than the continuum.
 We also disprove  (in ZFC) a conjecture of Hurewicz
 which is analogous to  the Borel conjecture.  Finally,
 we answer several questions from Part I concerning partition
 properties of covers.
 \footnote{
 AMS Classification: 03E05 04A20 54D20

 $\;$ Key words: Rothberger property $C^{\prime\prime}$,
 Gerlits-Nagy property
 $\gamma$-sets,  Hurewicz property,
 Menger property, $\gamma$--cover, $\omega$--cover,
 Sierpi\'nski set, Lusin set.
           }
\end{abstract}

\begin{center}{\bf Introduction}
\end{center}

Many topological properties of spaces have been defined or characterized
in terms of the properties of open coverings of these spaces.
Popular among such properties are the properties
introduced by Gerlits and Nagy \cite{G-N}, Hurewicz
\cite{Hu}, Menger \cite{Me} and Rothberger \cite{Ro}.
These are all defined in terms of the possibility of extracting
from a given sequence of open covers of some sort, an open cover of some
(possibly   different) sort.

In Scheepers \cite{S} it was shown that
when one systematically studies the definitions involved and inquires whether
other natural variations of the defining procedures produce any new classes
of sets which have mathematically interesting properties, an aesthetically
pleasing picture emerges.  In \cite{S} the basic implications were
established.  It was left open whether these were the only implications.

  Let $X$ be a topological space. By a ``cover'' for $X$ we always
mean ``countable open cover''. Since we are primarily interested in
separable metrizable (and hence Lindel\"of) spaces, the restriction
to countable covers does not lead to a loss of generality.
  A cover ${\cal U}$ of $X$ is said to be
\begin{enumerate}
 \item{{\em large $\lm$} if for each $x$ in $X$ the set
      $\{U\in {\cal U}:x\in U\}$ is infinite;}
 \item{{\em an $\omega$--cover $\om$} if $X$ is not in ${\cal U}$ and for
       each finite subset $F$ of $X$, there is
       a set $U\in{\cal U}$ such that $F\subset U$;}
 \item{{\em a $\gamma$--cover $\ga$} if it is infinite and for each $x$ in
       $X$ the set $\{U\in {\cal U}:x\not\in U\}$ is finite.}
\end{enumerate}

   We shall use the symbols $\op$, $\lm$, $\om$ and
$\ga$ to denote the collections of all open, large, $\omega$ and 
$\gamma$--covers respectively, of $X$. Let ${\cal A}$ and ${\cal B}$
each be one of these four classes.  We consider the
following three
``procedures'', $\sone$, $\sfin$ and $\ufin$, for obtaining covers in
${\cal B}$ from covers in ${\cal A}$:
\begin{enumerate}
\item{$\sone({\cal A},{\cal B})$: For a sequence
     $({\cal U}_n:n=1,2,3,\dots)$ of
     elements of ${\cal A}$, select for each $n$ a set
     $U_n\in{\cal U}_n$ such that
     $\{U_n:n=1,2,3,\dots\}$ is a member of ${\cal B}$;}
\item{$\sfin({\cal A},{\cal B})$: For a sequence
     $({\cal U}_n:n=1,2,3,\dots)$ of
     elements of ${\cal A}$, select for each $n$ a finite set
     ${\cal V}_n\subset{\cal U}_n$ such that
     $\cup_{n=1}^{\infty}{\cal V}_n$ is an element of ${\cal B}$;}
\item{$\ufin({\cal A},{\cal B})$: For a sequence
     $({\cal U}_n:n=1,2,3,\dots)$ of elements of ${\cal A}$,
     select for each $n$ a finite set ${\cal V}_n\subset{\cal U}_n$
     such that
     $\{\cup{\cal V}_n:n=1,2,3,\dots\}$ is a member of ${\cal B}$
     or\footnote{This is similar to the  $*$ convention of \cite{S}.}
     there exists an $n$ such that $\cup{\cal V}_n=X$.}
\end{enumerate}

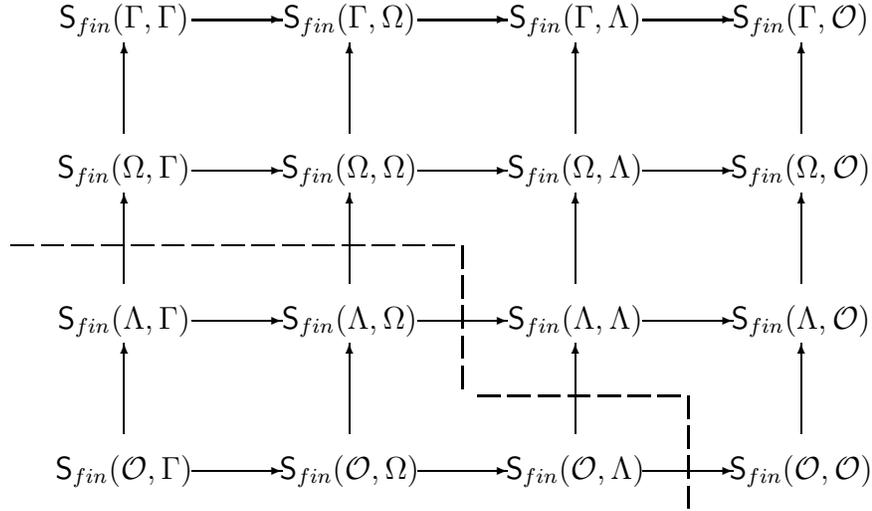
\begin{figure}
\unitlength=1.00mm
\begin{picture}(80.00,80.00)(-12,30)
\multiput(19,40)(0,20){4}{\multiput(0,0)(30,0){3}{\vector(1,0){12}}}
\multiput(10,45)(30,0){4}{\multiput(0,0)(0,20){3}{\vector(0,1){12}}}
\multiput(-5,70)(4,0){15}{\line(1,0){3}}
\multiput(57,50)(4,0){7}{\line(1,0){3}}
\multiput(55,70)(0,-4){5}{\line(0,-1){3}}
\multiput(85,50)(0,-4){4}{\line(0,-1){3}}
\put(10.00,40.00){\makebox(0,0)[cc]{$\sfin(\op,\ga)$}}
\put(40.00,40.00){\makebox(0,0)[cc]{$\sfin(\op,\om)$}}
\put(70.00,40.00){\makebox(0,0)[cc]{$\sfin(\op,\lm)$}}
\put(100.00,40.00){\makebox(0,0)[cc]{$\sfin(\op,\op)$}}
\put(10.00,60.00){\makebox(0,0)[cc]{$\sfin(\lm,\ga)$}}
\put(40.00,60.00){\makebox(0,0)[cc]{$\sfin(\lm,\om)$}}
\put(70.00,60.00){\makebox(0,0)[cc]{$\sfin(\lm,\lm)$}}
\put(100.00,60.00){\makebox(0,0)[cc]{$\sfin(\lm,\op)$}}
\put(10.00,80.00){\makebox(0,0)[cc]{$\sfin(\om,\ga)$}}
\put(40.00,80.00){\makebox(0,0)[cc]{$\sfin(\om,\om)$}}
\put(70.00,80.00){\makebox(0,0)[cc]{$\sfin(\om,\lm)$}}
\put(100.00,80.00){\makebox(0,0)[cc]{$\sfin(\om,\op)$}}
\put(10.00,100.00){\makebox(0,0)[cc]{$\sfin(\ga,\ga)$}}
\put(40.00,100.00){\makebox(0,0)[cc]{$\sfin(\ga,\om)$}}
\put(70.00,100.00){\makebox(0,0)[cc]{$\sfin(\ga,\lm)$}}
\put(100.00,100.00){\makebox(0,0)[cc]{$\sfin(\ga,\op)$}}
\end{picture}
\caption{Basic diagram for $\sfin$ \label{sfinfig}}
\end{figure}

For ${\sf G}$ one of these three procedures, let us say that
a space has property
${\sf G}({\cal A},{\cal B})$ if for every sequence of elements
of ${\cal A}$, one can
obtain an element of ${\cal B}$ by means of procedure ${\sf G}$.
Letting ${\cal A}$ and
${\cal B}$ range over the set $\{\op, \lm, \om, \ga \}$, we see that
for each ${\sf G}$ there are potentially sixteen classes of
spaces of the form ${\sf G}({\cal A},{\cal B})$.
Each of our properties is monotone decreasing in the first coordinate
and increasing in the second, hence we get the following diagram
(see figure \ref{sfinfig})
for ${\sf G}=\sfin$.

It also
is easily checked that $\sfin(\lm,\om)$ and
$\sfin(\op,\lm)$ are impossible for nontrivial
$X$.  Hence the five classes in the lower left corner
are eliminated.  The same follows for the stronger
property $\sone$.  In the case of $\ufin$ note that
$\ufin(\op,\cdot)$ is equivalent to $\ufin(\ga,\cdot)$
because given an open cover $\{U_n: n\in\omega\}$ we may replace
it by the $\gamma$-cover, $\{\cup_{i<n} U_i: n\in\omega\}$.
This means the diagram of $\ufin$ reduces to any of its rows.
Now clearly $\sone$ implies $\sfin$. Also it is clear that
$\sfin(\ga,{\cal A})\rightarrow\ufin(\ga,{\cal A})$
for ${\cal A}=\ga,\om,\op$.
The implication $$\sfin(\ga,\lm)\rightarrow\ufin(\ga,\lm)$$
is also true, but takes a little thought since when we take
finite unions we might not get distinct sets.  To prove it,
assume ${\cal U}_n$ are $\gamma$--covers of $X$ with no
finite subcover.  Applying $\sfin$ we get
a sequence of finite ${\cal V}_n\subseteq {\cal U}_n$ such that
for any $x$ there exists infinitely many $n$ such that
$x\in \cup {\cal V}_n$.
But since the ${\cal U}_m$'s have no finite subcover
we can inductively choose finite ${\cal W}_n$ with
$${\cal V}_n\subseteq {\cal W}_n\subseteq {\cal U}_n$$
and
$\cup{\cal W}_n\not=\cup{\cal W}_m$ for any $m<n$. Hence
$\bigcup\{{\cal W}_n:n=1,2,3,\ldots\}$ is a large cover of $X$.

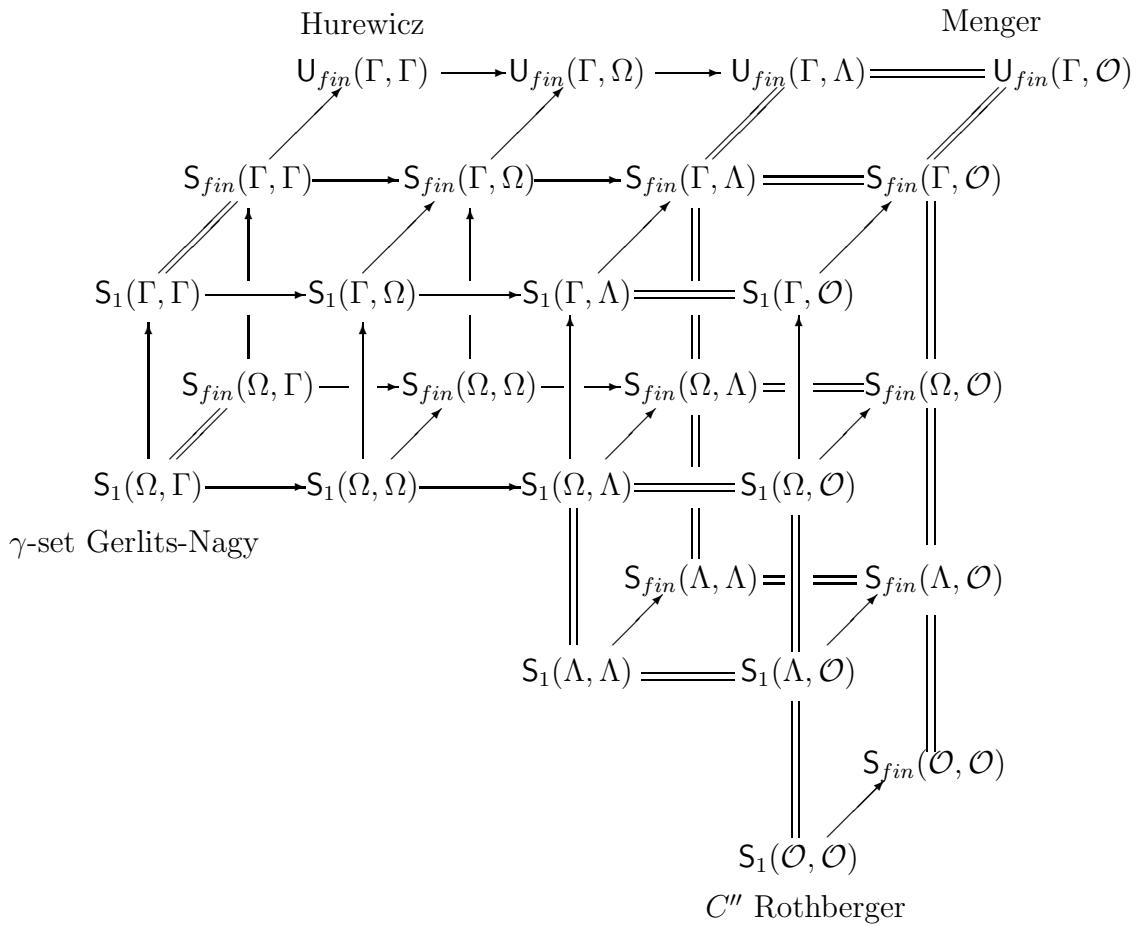
\begin{figure}
\unitlength=.95mm
\begin{picture}(141.00,145.00)(0,0)
\put(104.00,20.00){\makebox(0,0)[cc]{$\sone(\op,\op)$}}
\put(104.00,46.00){\makebox(0,0)[cc]{$\sone(\lm,\op)$}}
\put(73.00,46.00){\makebox(0,0)[cc]{$\sone(\lm,\lm)$}}
\put(104.00,72.00){\makebox(0,0)[cc]{$\sone(\om,\op)$}}
\put(73.00,72.00){\makebox(0,0)[cc]{$\sone(\om,\lm)$}}
\put(43.00,72.00){\makebox(0,0)[cc]{$\sone(\om,\om)$}}
\put(13.00,72.00){\makebox(0,0)[cc]{$\sone(\om,\ga)$}}
\put(104.00,99.00){\makebox(0,0)[cc]{$\sone(\ga,\op)$}}
\put(73.00,99.00){\makebox(0,0)[cc]{$\sone(\ga,\lm)$}}
\put(43.00,99.00){\makebox(0,0)[cc]{$\sone(\ga,\om)$}}
\put(13.00,99.00){\makebox(0,0)[cc]{$\sone(\ga,\ga)$}}

\put(123.00,33.00){\makebox(0,0)[cc]{$\sfin(\op,\op)$}}
\put(123.00,59.00){\makebox(0,0)[cc]{$\sfin(\lm,\op)$}}
\put(89.00,59.00){\makebox(0,0)[cc]{$\sfin(\lm,\lm)$}}
\put(123.00,86.00){\makebox(0,0)[cc]{$\sfin(\om,\op)$}}
\put(89.00,86.00){\makebox(0,0)[cc]{$\sfin(\om,\lm)$}}
\put(58.00,86.00){\makebox(0,0)[cc]{$\sfin(\om,\om)$}}
\put(27.00,86.00){\makebox(0,0)[cc]{$\sfin(\om,\ga)$}}
\put(27.00,115.00){\makebox(0,0)[cc]{$\sfin(\ga,\ga)$}}
\put(58.00,115.00){\makebox(0,0)[cc]{$\sfin(\ga,\om)$}}
\put(89.00,115.00){\makebox(0,0)[cc]{$\sfin(\ga,\lm)$}}
\put(123.00,115.00){\makebox(0,0)[cc]{$\sfin(\ga,\op)$}}

\put(141.00,130.00){\makebox(0,0)[cc]{$\ufin(\ga,\op)$}}
\put(43.00,130.00){\makebox(0,0)[cc]{$\ufin(\ga,\ga)$}}
\put(73.00,130.00){\makebox(0,0)[cc]{$\ufin(\ga,\om)$}}
\put(104.00,130.00){\makebox(0,0)[cc]{$\ufin(\ga,\lm)$}}

\put(105.00,13.00){\makebox(0,0)[cc]{$C^{\prime\prime}$
Rothberger}}
\put(11.00,64.00){\makebox(0,0)[cc]{$\gamma$-set Gerlits-Nagy}}
\put(131.00,137.00){\makebox(0,0)[cc]{Menger}}
\put(43.00,137.00){\makebox(0,0)[cc]{Hurewicz}}

\put(054.00,130.00){\vector(1,0){9.00}}
\put(084.00,130.00){\vector(1,0){9.00}}
\put(114.00,130.50){\line(1,0){16.00}}
\put(114.00,129.50){\line(1,0){16.00}}

\put(061.00,118.00){\vector(1,1){10.00}}
\put(092.00,118.00){\line(1,1){10.00}}
\put(091.00,118.00){\line(1,1){10.00}}
\put(030.00,118.00){\vector(1,1){10.00}}
\put(122.00,118.00){\line(1,1){10.00}}
\put(123.00,118.00){\line(1,1){10.00}}

\put(036.00,115.00){\vector(1,0){12.00}}
\put(067.00,115.00){\vector(1,0){12.00}}
\put(099.00,115.50){\line(1,0){14.00}}
\put(099.00,114.50){\line(1,0){14.00}}

\put(014.50,102.00){\line(1,1){10.00}}
\put(015.50,102.00){\line(1,1){10.00}}
\put(043.00,102.00){\vector(1,1){10.00}}
\put(076.00,102.00){\vector(1,1){10.00}}
\put(107.00,102.00){\vector(1,1){10.00}}

\put(090.00,101.00){\line(0,1){10.00}}
\put(089.00,101.00){\line(0,1){10.00}}
\put(058.00,101.00){\vector(0,1){10.00}}
\put(027.00,101.00){\vector(0,1){10.00}}

\put(021.00,099.00){\vector(1,0){14.00}}
\put(051.00,099.00){\vector(1,0){14.00}}
\put(081.00,099.50){\line(1,0){14.00}}
\put(081.00,098.50){\line(1,0){14.00}}

\put(027.00,090.00){\line(0,1){6.00}}
\put(058.00,090.00){\line(0,1){6.00}}
\put(090.00,090.00){\line(0,1){6.00}}
\put(089.00,090.00){\line(0,1){6.00}}
\put(123.00,090.00){\line(0,1){22.00}}
\put(122.00,090.00){\line(0,1){22.00}}

\put(037.00,086.00){\line(1,0){4.00}}
\put(045.00,086.00){\vector(1,0){3.00}}
\put(068.00,086.00){\line(1,0){3.00}}
\put(074.00,086.00){\vector(1,0){5.00}}
\put(099.00,086.50){\line(1,0){3.00}}
\put(099.00,085.50){\line(1,0){3.00}}
\put(106.00,086.50){\line(1,0){7.00}}
\put(106.00,085.50){\line(1,0){7.00}}

\put(016.00,076.00){\line(1,1){7.00}}
\put(017.00,076.00){\line(1,1){7.00}}
\put(047.00,076.00){\vector(1,1){7.00}}
\put(077.00,076.00){\vector(1,1){7.00}}
\put(107.00,076.00){\vector(1,1){7.00}}

\put(013.00,076.00){\vector(0,1){19.00}}
\put(043.00,076.00){\vector(0,1){19.00}}
\put(072.00,076.00){\vector(0,1){20.00}}
\put(090.00,075.00){\line(0,1){7.00}}
\put(089.00,075.00){\line(0,1){7.00}}
\put(104.00,076.00){\vector(0,1){20.00}}

\put(021.00,072.00){\vector(1,0){14.00}}
\put(051.00,072.00){\vector(1,0){14.00}}
\put(081.00,072.50){\line(1,0){14.00}}
\put(081.00,071.50){\line(1,0){14.00}}

\put(089.00,062.00){\line(0,1){7.00}}
\put(090.00,062.00){\line(0,1){7.00}}
\put(122.00,064.00){\line(0,1){19.00}}
\put(123.00,064.00){\line(0,1){19.00}}

\put(099.00,059.50){\line(1,0){3.00}}
\put(099.00,058.50){\line(1,0){3.00}}
\put(106.00,059.50){\line(1,0){6.00}}
\put(106.00,058.50){\line(1,0){6.00}}

\put(078.00,050.00){\vector(1,1){7.00}}
\put(108.00,050.00){\vector(1,1){7.00}}

\put(073.00,050.00){\line(0,1){19.00}}
\put(072.00,050.00){\line(0,1){19.00}}
\put(104.00,049.00){\line(0,1){19.00}}
\put(103.00,049.00){\line(0,1){19.00}}

\put(082.00,046.00){\line(1,0){13.00}}
\put(082.00,045.00){\line(1,0){13.00}}

\put(123.00,035.00){\line(0,1){19.00}}
\put(122.00,035.00){\line(0,1){19.00}}

\put(104.00,023.00){\line(0,1){19.00}}
\put(103.00,023.00){\line(0,1){19.00}}

\put(108.00,023.00){\vector(1,1){8.00}}

\end{picture}
\caption{ Full 3d diagram \label{3dim}}
\end{figure}

\begin{figure}
\unitlength=.95mm
\begin{picture}(140.00,100.00)(10,10)
\put(20.00,20.00){\makebox(0,0)[cc]
{\shortstack {$\sone(\om,\ga)$\\ $\{\}$         } }}
\put(60.00,20.00){\makebox(0,0)[cc]
{\shortstack {$\sone(\om,\om)$\\ $\{ H   \} $   } }}
\put(100.00,20.00){\makebox(0,0)[cc]
{\shortstack {$\sone(\op,\op)$\\ $\{L,H   \}$   } }}
\put(20.00,60.00){\makebox(0,0)[cc]
{\shortstack {$\sone(\ga,\ga)$\\ $\{ S   \} $   } }}
\put(60.00,60.00){\makebox(0,0)[cc]
{\shortstack {$\sone(\ga,\om)$\\ $\{ H,S  \}$   } }}
\put(100.00,60.00){\makebox(0,0)[cc]
{\shortstack {$\sone(\ga,\lm)$\\ $\{L,H,S  \}$  } }}
\put(80.00,40.00){\makebox(0,0)[cc]
{\shortstack {$\sfin(\om,\om)$\\ $\{ C,H  \}$   } }}
\put(80.00,80.00){\makebox(0,0)[cc]
{\shortstack {$\sfin(\ga,\om)$\\ $\{ C,S,H \}$  } }}
\put(60.00,100.00){\makebox(0,0)[cc]
{\shortstack {$\ufin(\ga,\ga)$\\ $\{C,S   \}$   } }}
\put(100.00,100.00){\makebox(0,0)[cc]
{\shortstack {$\ufin(\ga,\om)$\\ $\{C,S,H  \}$  } }}
\put(140.00,100.00){\makebox(0,0)[cc]
{\shortstack {$\ufin(\ga,\op)$\\ $\{C,S,L,H \}$ } }}

\put(70.00,100.00){\vector(1,0){20.00}}
\put(110.00,100.00){\vector(1,0){20.00}}
\put(86.00,85.00){\vector(1,1){10.00}}
\put(25.00,65.00){\vector(1,1){30.00}}
\put(105.00,65.00){\vector(1,1){30.00}}
\put(65.00,65.00){\vector(1,1){10.00}}
\put(28.00,61.00){\vector(1,0){20.00}}
\put(70.00,61.00){\vector(1,0){20.00}}
\put(80.00,45.00){\line(0,1){12.00}}
\put(80.00,64.00){\vector(0,1){11.00}}
\put(20.00,25.00){\vector(0,1){29.00}}
\put(60.00,25.00){\vector(0,1){29.00}}
\put(64.00,27.00){\vector(1,1){8.00}}
\put(100.00,25.00){\vector(0,1){29.00}}
\put(28.00,20.00){\vector(1,0){20.00}}
\put(70.00,20.00){\vector(1,0){20.00}}
\end{picture}
\caption{After removing equivalent classes \label{cshl}}
\end{figure}
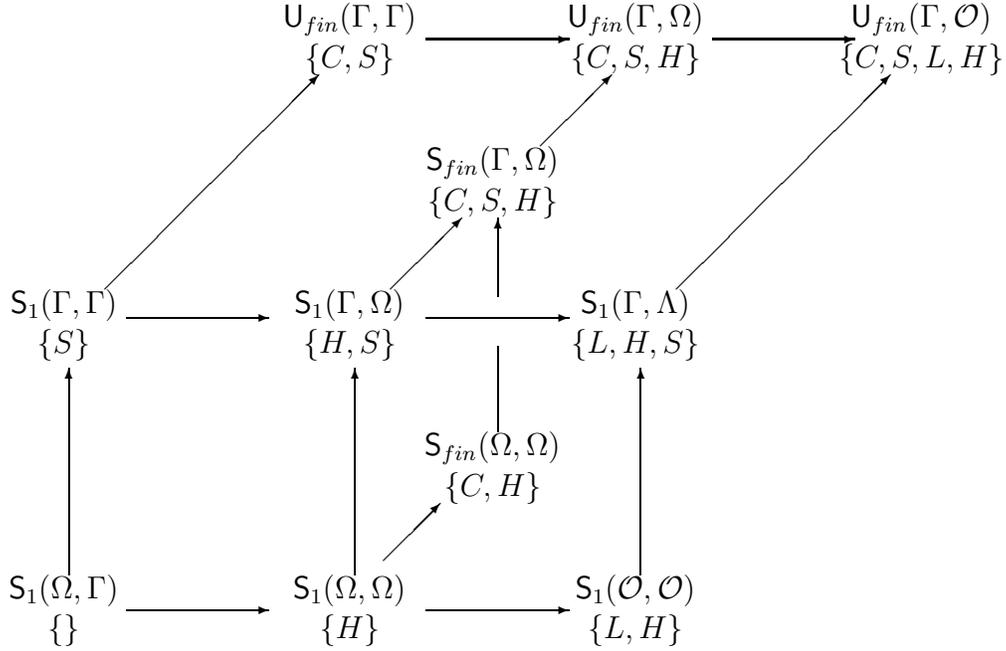

In the three dimensional diagram of figure~\ref{3dim}
the double lines indicate that the two properties are equivalent.
The proof of these equivalences can be found in either
Scheepers \cite{S} or section~\ref{equiv}  of this paper.
After removing duplications we obtain figure \ref{cshl}.

For this diagram, we have provided four examples $\{C,S,H,L\}$
which show that practically no other implications can hold.
$C$ is the Cantor set ($2^\omega$), $S$ is a special Sierpi\'nski set
such the $S+S$ can be mapped continuously onto the irrationals,
$L$ is a special
Lusin set such that $L+L$ can be mapped continuously onto the irrationals,
and $H$ is a generic Lusin set.  Thus the only remaining
problems are:

\begin{problem} Is $\ufin(\ga,\om)=\sfin(\ga,\om)$?
\end{problem}

\begin{problem} And if not, does $\ufin(\ga,\ga)$ imply
$\sfin(\ga,\om)$?
\end{problem}

All of our examples are subsets of the real line, but only one of
them (the Cantor set) is a ZFC example. Thus, the following problem
arises:
\begin{problem} Are there ZFC examples of (Lindel\"of) topological
spaces which show that none of the arrows in figure \ref{cshl} can
even be consistently reversed?
\end{problem}

The paper is organized as follows:

In section \ref{equiv} we prove the equivalences of our
properties indicated in figure \ref{3dim}.
We prove that $\sone(\ga,\ga)=\sfin(\ga,\ga)$ and
$\sfin(\lm,\lm)=\sfin(\ga,\lm)$.
The other equivalences are either trivial or were proved in
Scheepers \cite{S}. In section~\ref{examp}
we present the four examples C,S,H,L indicated in figure \ref{cshl}.
In section~\ref{pres} we study the preservation of these
families under the taking of finite powers and other operations.

In section~\ref{card} we study for each of these eleven
families the cardinal
$$\non({\cal X})=\{\min |X|: X\notin {\cal X}\}.$$
We show that each is equal to either $\goth b$, $\goth d$,
$\goth p$, or the covering number of the meager ideal $\covmeag$.
We also show that ${\goth r}=\non(\split(\lm,\lm))$ and
${\goth u}=\non(\split(\om,\om))$. ($\split({\cal A},{\cal B})$
holds iff every infinite cover from ${\cal A}$ can be split
into two covers from ${\cal B}$).

In section \ref{hur} we give a ZFC counterexample to a conjecture
of Hurewicz by showing that there exists an uncountable set
of reals in $\ufin(\ga,\ga)$ which is not $\sigma$-compact.
We also show the any $\ufin(\ga,\ga)$ set which does not
contain a perfect set is perfectly meager.

In section \ref{ramsey} we  consider other properties
from Scheepers \cite{S} and settle some
questions about Ramsey-like properties of covers that were left
open in \cite{S}.  We show that
$\sone(\om,\om)$ implies $Q(\om,\om)$ and hence
$$\sone(\om,\om) = {\sf P}(\om,\om)+{\sf Q}(\om,\om).$$
We also show that $\om\to \lceil\om\rceil^2_2$
is equivalent to $\sfin(\om,\om)$.

\section{Equivalences}\label{equiv}

In this section we show  $\sone(\ga,\ga) = \sfin(\ga,\ga)$
and $\sfin(\ga,\lm)=\sfin(\lm,\lm)$.

The equivalence
$\sone(\om,\ga)$ with the $\gamma$-set property
(every $\omega$-cover contains a $\gamma$-cover) was shown
by Gerlits and Nagy \cite{G-N}.  But it is easy to see that
$\sfin(\om,\ga)$ implies the $\gamma$-set property and hence
$\sone(\om,\ga)=\sfin(\om,\ga)$.
In Scheepers \cite{S} (Cor 6) it was shown that
$\sone(\ga,\lm)=\sone(\ga,\op)$.

All of the other equivalences are either to the
Rothberger property $C^{\prime\prime}$ or
the Menger property.
For the Menger property, in Scheepers \cite{S} (Cor 5) it was shown that
$\sfin(\ga,\lm)=\ufin(\ga,\op)$.  $\sfin(\ga,\lm)=\sfin(\lm,\lm)$
by Theorem \ref{mengerst1}
and note also that $\sfin(\op,\op)$ easily follows from
$\ufin(\ga,\op)$ and hence all nine classes
(see figure \ref{3dim}) are equivalent to the Menger property.
In \cite{S} (Thm 17) it was shown that
all five classes (see figure \ref{3dim})
are equivalent to the Rothberger property ($C^{\prime\prime}$).

\begin{theorem}\label{sone_eq_sfin} $\sone(\ga,\ga) = \sfin(\ga,\ga)$.
\end{theorem}
\pf

Note that
the class $\sone(\ga,\ga)$ is contained 
in the class $\sfin(\ga,\ga)$. The difficulty with
showing that these two classes 
are in fact equal is as follows: when we are allowed to choose
finitely many elements per
$\gamma$--cover, we are allowed to also pick no elements; for
$\sone(\ga,\ga)$ we are required to
choose an element per $\gamma$--cover.

Let $X$ be a space which has property
$\sfin(\ga,\ga)$, and for
  each $n$ let ${\cal U}_n$ be a $\gamma$--cover of $X$,
enumerated bijectively
   as  $(U^n_1,U^n_2,U^n_3,...)$.

  For each $n$ define ${\cal V}_n$ to be
$\{V^n_1,V^n_2,V^n_3,...\}$, where

\[V^n_k= U^1_k\cap U^2_k\cap U^3_k\cap \dots\cap U^n_k.
\]
  For each $n$ , ${\cal V}_n$ is a $\gamma$--cover: For
  fix $n$.  For each $x$, and for each $i\in\{1,\dots,n\}$ there
exists an $N_i$
  such that $x$ is in $U^i_m$ for all $m>N_i$.  But then $x$ is in
$V^n_m$ for all
  $m> \max\{N_i:i=1,\dots,n\}$.

  Now apply $\sfin(\ga,\ga)$ to $({\cal V}_n :
n=1,2,\dots)$. We get a   sequence
             $$({\cal W}_n:n\in\omega)$$
such that ${\cal W}_n$ is a finite subset
of ${\cal V}_n$ for each $n$, and such that
$\cup_{n=1}^{\infty}{\cal W}_n$ is
a $\gamma$--cover of $X$.

  Choose a an increasing sequence $n_1<n_2<\dots<n_k<\dots$ such
that for each $j$,
  ${\cal W}_{n_j}\setminus\cup_{i<j}{\cal W}_{n_i}$ is nonempty.
This is possible because each ${\cal W}_n$ is finite, while
  the union of these sets, being a $\gamma$--cover of $X$, is
infinite.
  For each $j$, choose $m_j$ such that $V^{n_j}_{m_j}$ is an
element of ${\cal W}_{n_j}\setminus\cup_{i<j}{\cal W}_{n_i}$.
  Then $\{V^{n_k}_{m_k}: k=1,2,\dots\}$ is a $\gamma$--cover of
$X$.

  For each $n$ in $(n_k,n_{k+1}]$ we define $p_n=m_{k+1}$. Then
$\{U^n_{p_n}: n=1,2,\dots\}$ is a $\gamma$ cover  of $X$.
\qed

\begin{theorem}\label{mengerst1}
$\sfin(\ga,\lm)=\sfin(\lm,\lm)$.
\end{theorem}
\pf  Since $\ga\subseteq\lm$, it is clear that
$\sfin(\lm,\lm)$ implies
$\sfin(\ga,\lm)$.  We prove the other implication.

Assume  $\sfin(\ga,\lm)$
and let $({\cal U}_n:n\in\omega)$ be a sequence of large
covers of $X$.  Without loss of generality we may assume
that for every finite $F\subseteq \bigcup_{n\in\omega}{\cal U}_n$
that ${\cal U}_k\cap F=\emptyset$ for all but finitely many $k$.
(This can be accomplished by throwing out finitely many elements
from each ${\cal U}_n$.)

For each $n$ enumerate ${\cal U}_n$
bijectively as $(U^n_k:k\in\omega)$, and
define
$$V^n_m =\bigcup\{U^n_i: i<m\}.$$

Since each ${\cal V}_n=(V^n_m:m\in\omega)$ is a nondecreasing
open cover of $X$, either there exists $m_n$ such that
$V_{m_n}^n=X$ or
${\cal V}_n$ is a $\gamma$-cover.  So there must be an infinite
$A$ for which one or the other always occurs.  Suppose
${\cal V}_n$ is a $\gamma$-cover
for every $n\in A$ .
Apply  $\sfin(\ga,\lm)$ to obtain
${\cal W}_n$ a finite subset of ${\cal V}_n$
such that $\bigcup\{{\cal W}_n:n\in A\}$  is a large cover of $X$.
Let ${\cal P}_n$ be a finite subset of ${\cal U}_n$ such
that every element of ${\cal W}_n$ is a union of elements
of ${\cal P}_n$.  Since ${\cal P}_n$ is disjoint from all
but finitely many of the ${\cal U}_k$, it follows that
$\bigcup\{{\cal P}_n: n \in A\}$ is a large cover of
$X$.  In the case that $V_{m_n}^n=X$ for every $n\in A$ just
take ${\cal P}_n=\{U^n_i: i<m_n\}$ and the same argument works.
\qed

\newpage

\section{Examples}\label{examp}

\subsec {The Cantor set C}

It is easy to check that every $\sigma$--compact space
(the union of countably many compact sets)
belongs to $\ufin(\ga,\ga)$ and $\sfin(\om,\om)$.
We also show that the Cantor set, $2^\omega$, is
not in the class $\sone(\ga,\lm)$.

For the sake of conciseness, let us introduce the following
notion. An open cover ${\cal U}$ of a topological space $X$ is a
{\em $k$--cover} iff
there is for every $k$--element subset of $X$ an element of
${\cal U}$  which covers that set.

\begin{lemma}\label{cptomeg} Let $k$ be a positive integer. Every
$\omega$--cover of a compact space contains a finite subset
which is a $k$--cover for the space.
\end{lemma}
\pf
   Let ${\cal U}$ be an $\omega$--cover of the compact space $X$
and let $k$
be a positive integer. Then the set ${\cal V} = \{U^k:U\in{\cal U}\}$ of
$k$--th powers of elements of ${\cal U}$ is a collection of open subsets
of $X^k$, and it is a cover of $X^k$ since ${\cal U}$ is an
$\omega$--cover of $X$. Since $X$ is compact, so is $X^k$. Thus
there is a finite subset of ${\cal V}$ which covers $X^k$, say
$\{U^k_1,\dots,U^k_n\}$. But then $\{U_1,\dots,U_n\}$ is a
$k$--cover of$X$.
\qed

\begin{theorem} Every $\sigma$--compact topological space is a
member of both the class ${\sfin}(\om,\om)$ and $\ufin(\ga,\ga)$.
\end{theorem}
\pf Let $X$ be a $\sigma$--compact space, and
write
$X=\cup_{n\in\omega}K_n$ where
$$K_0\subseteq K_1\subseteq \dots \subseteq K_n\subseteq$$
is a sequence of compact subsets of $X$.

Let $({\cal U}_n:n\in\omega)$
be a sequence of $\omega$--covers of $X$.
For each
$n$ apply Lemma \ref{cptomeg} to the $\omega$--cover ${\cal U}_n$ of the
space $K_n$, to find a finite subset ${\cal V}_n$ of ${\cal U}_n$ which
is an $n$--cover of $K_n$. Then $\cup_{n\in\omega}{\cal V}_n$
is an  $\omega$--cover of $X$.  This shows that $X$ has property
$\sfin(\om,\om)$.

Now suppose that $({\cal U}_n:n\in\omega)$ is a sequence of
$\gamma$-covers of $X$.  Since any infinite subset of a
$\gamma$-cover is a $\gamma$-cover we may assume that our
covers are disjoint.  Since each $K_n$ is compact we may choose
${\cal V}_n\in [{\cal U}_n]^{<\omega}$ so that
$K_n\subseteq \cup{\cal V}_n$.
Either there exists
$n$ such that $\cup{\cal V}_n=X$ or
$\{\cup{\cal V}_n:n\in\omega\}$ is infinite and hence
a $\gamma$-cover of $X$.   It follows that
$X$ has property $\ufin(\ga,\ga)$.
\qed

\begin{theorem}
  The Cantor set, $C=2^\omega$, is not in the class
  $\sone(\ga,\lm)$.
\end{theorem}

\pf
There exists an $\omega\times\omega$--matrix
$(A^m_n:m,n<\omega)$ of
closed subsets of the Cantor set such that
\begin{enumerate}
\item for each fixed $m\in\omega$ the sets $A^m_n$ for $n<\omega$ are
 pairwise disjoint and \label{disc1}
\item whenever $m_1<m_2<\dots<m_k$ and $n_1,n_2,\dots,n_k$ are \label{disc2}
given, then  $A^{m_1}_{n_1}\cap\dots\cap A^{m_k}_{n_k}\neq\emptyset$.
\end{enumerate}
To see that such a matrix exists
think of the Cantor set as the homeomorphic space
$2^{\omega\times\omega}$
 instead.

   Let $\la x_n,\ n<\omega\ra$ be a sequence of
pairwise distinct elements of
   $2^{\omega}$. Also, for each $m$, let
\[\pi_m:2^{\omega\times\omega}\rightarrow \mbox{$2^{\omega}$}
\]
   be defined so that for each $y$ in $2^{\omega\times\omega}$
   and for each $m$,
$\pi_m(y)(n) = y(m,n)$.
   Then $\pi_m$ is continuous. We now define our matrix.

   For each $m$ and $n$ we define
\[A^m_n = \{y\in 2^{\omega\times\omega}:\pi_m(y) = x_n\}
\]
   Each row of the matrix is pairwise disjoint since the $x_n$'s
are pairwise distinct.
   Each entry of the matrix is a closed set since each $\pi_m$ is
continuous. We must
   still verify property 2. Thus, let
$(m_1,n_1),\dots,(m_k,n_k)$ be given
   such that $m_1<\dots<m_k$. Then the element $y$ of
$2^{\omega\times\omega}$ which is
defined so that for each~$i$
\[\mbox{for each $j$, }y(m_i,j)=x_{n_i}(j)
\]
is a member of the set $A^{m_1}_{n_1}\cap\dots\cap
A^{m_k}_{n_k}$, whence this
intersection is nonempty.

 For each  $m$ put
${\cal U}_m = \{2^\omega\setminus A^m_n:n<\omega\}$. Then by
property \ref{disc1}
we see that each ${\cal U}_m$ is a $\gamma$--cover of
$2^\omega$.

   For each $m$ choose a $U^m_{n_m}$ from ${\cal U}_m$. Then
$U^m_{n_m} =
 2^\omega\setminus A^m_{n_m}$. By the property \ref{disc2}
and  the fact that all the $A^m_n$'s are
compact,
we see that the intersection $\cap_{m<\omega}A^m_{n_m}$ is
nonempty. But then
not only is $\{U^m_{n_m}:m<\omega\}$ not a large cover of
$2^{\omega\times\omega}$,
it is not even a cover of $2^{\omega\times\omega}$.
\qed
It follows from these two theorems that the Cantor set C must
lie exactly in those classes indicated in figure \ref{cshl} in
our introduction.

\begin{theorem} No uncountable $F_{\sigma}$ set of reals is
  in $\sone(\ga,\ga)$. \label{thm17}
\end{theorem}
\pf Such a set contains an uncountable compact perfect set. The
Cantor set is a continuous image of such perfect sets.
\qed

\subsec{The special Lusin set L}

Recall that a set $L$ of real numbers  is
said to be a {\em Lusin} set iff it is uncountable
but its intersection with every first category set of real
numbers is countable.
Sierpi\'nski \cite{sier} showed that assuming CH there
exists a Lusin set $L$ such that $L+L$ is the irrationals
(see also Miller \cite{survey} Thm 8.5).

We will construct similarly a Lusin set $L\subseteq \zz$ with
the property that $L+L=\zz$.  Here $\zz$ is the infinite product of the
ring of integers and addition is the usual pointwise addition,
i.e, $(x+y)(n)=x(n)+y(n)$.
Our construction is based on the
following simple fact:

\begin{lemma}\label{comeagersums} If $X$ is a comeager subset of
$\zz$, then for every
$x\in\zz$ there are elements $a$ and $b$ of $X$ such that
$a+b=x$.
\end{lemma}
\pf
Since multiplication
by $-1$ and translation by $x$ are homeomorphisms, the
 set
 $$x-X = \{x-y:y\in X\}$$
is also comeager. But then $X\cap (x-X)$ is non--empty.
Let $z$ be an element of this intersection. Then $z=a$ for some
$a$ in $X$, and $z=x-b$ for some $b$ in $X$. The Lemma follows.
\qed

\begin{lemma} [CH]\label{lusinset} There is a
Lusin set $L\subseteq\zz$
such that $L+L=\zz$.
\end{lemma}
\pf
Let
$(M_{\alpha}:\alpha<\omega_1)$ bijectively list all first
category $F_{\sigma}$--subsets of $\zz$. Let
$(r_{\alpha}:\alpha<\omega_1)$
bijectively list $\zz$. Using Lemma \ref{comeagersums}, choose elements
$x_{\alpha}, y_{\alpha}$ from $\zz$ subject to the following rules:
\begin{enumerate}
\item For each $\alpha$, $r_{\alpha} = x_{\alpha} + y_{\alpha}$, and
\item $x_{\alpha}$ and $y_{\alpha}$ are not elements of
$\cup_{\beta\leq \alpha}M_{\beta}\cup\{x_{\beta},y_\beta:\beta<\alpha\}$
\end{enumerate}
  Letting $L$ be the set
$\{x_{\alpha}:\alpha<\omega_1\}\cup\{y_{\alpha}:\alpha<\omega_1\}$
completes the proof.
\qed

For a proof of the following result see Rothberger \cite{Ro}.

\begin{theorem}
  (Rothberger) Every Lusin set has property
  $\sone(\op,\op)=C^{\prime\prime}$.
\end{theorem}

\begin{theorem} \label{specluz}
If $L$ is our special Lusin set (i.e., $L+L=\zz$),
then $L$ does not satisfy $\ufin(\ga,\om)$.
\end{theorem}

\pf
Let $\{{\cal U}_n: n\in\omega\}$ be the sequence of open covers
defined by
$${\cal U}_n=\{U_{n,k}:{k\in \omega}\}$$
where
$$U_{n,k}=\{f\in \zz : |f(n)|\leq k\}.$$
Then each ${\cal U}_n $ is a
$\gamma$-cover of ${L}$. Let $\{{\cal V}_n: n\in\omega\}$ be
a sequence such that
${\cal V}_n\in [{\cal U}_n]^{<\omega}$, and let $h\in \oo$
be such that
$$h(n)>2\cdot\max\{ k:U_{n,k}\in {\cal V}_n\}$$
for all $n\in \omega$.
Let $f,g\in { L}$ be such that $h=f+g$. Then
$$\max\{|f(n)|,|g(n)|\}\geq {1\over 2} h(n)$$
for all $n\in \omega$, and hence $\{f,g\}\not\subseteq \cup{\cal V}_n$
for  any $n\in \omega$.
\qed

\newpage

\subsec{The special Sierpi\'nski set S}

   A Sierpi\'nski set is an uncountable subset of the real line which
has countable intersection with every set of Lebesgue measure zero.
In Theorem 7 of Fremlin and Miller \cite{F-M} it was shown that
every Sierpi\'nski set
belongs to the class $\ufin(\ga,\ga)$. It is well known that every
Borel image of a Sierpi\'nski set into the Baire space 
$^{\omega}\omega$ is bounded.
Sets with the property that every Borel image in the
Baire space is bounded were called $A_2$--sets in Bartoszynski
and Scheepers \cite{B-S}.

\begin{theorem} Every $A_2$--set (hence
every Sierpi\'nski set) belongs to $\sone(\ga,\ga)$.
\end{theorem}
\pf Let $X$ be an $A_2$--set, and let $({\cal U}_n:n\in\omega)$ be a
sequence of $\gamma$--covers of it. Enumerate each ${\cal U}_n$
 bijectively as $(U^n_m:m\in\omega)$.

Define a function $\Psi$ from $X$ to $\oo$ so that for each
$x\in X$ and for each $n$,
$$\Psi(x)(n) = \min\{m:(\forall k\geq m)(x\in U^n_k)\}.$$
Then $\Psi$ is a Borel function. Choose a strictly
increasing
function $g$ from $\oo$ which eventually dominates each
element of $\Psi[X]$. Then the sequence
$(U^n_{g(n)}:n\in\omega)$ is a $\gamma$--cover of $X$.
\qed

Clearly no Sierpi\'nski set is of measure zero and
since every $\sone(\op,\op)$ set is of measure zero, $X$ fails
to be $\sone(\op,\op)$. Therefore we have established the following theorem:

\begin{theorem}\label{sierp} If $X$ is a Sierpi\'nski set of reals,
then $X$ is $\sone(\ga,\ga)$ but not $\sone(\op,\op)$.
\end{theorem}

We call a Sierpi\'nski set $S$ special iff
$S+S$ is the set of irrationals.  (Here we are
using ordinary addition in the reals.)  Using an argument similar
to Lemma \ref{lusinset} one can show that assuming
CH there exists a special Sierpi\'nski set.

\begin{theorem}
  A special Sierpi\'nski set is not in the class $\sfin(\om,\om)$.
\end{theorem}

\pf
By Theorem \ref{contimage} all
our classes are closed under continuous images.  Note
that $S+S$ is the continuous image of $S\times S$.
Also  $\oo$ is not in $\ufin(\ga,\op)$ (see proof of
Theorem  \ref{specluz}).
Hence $\oo$ is not in $\sfin(\om,\om)$ and therefore $S\times S$
is not in $\sfin(\om,\om)$.
But by Theorem \ref{sfinomegprod} the class $\sfin(\om,\om)$ is
closed under finite products and therefore $S$
is not in the class $\sfin(\om,\om)$.
\qed

These results show that the special Sierpi\'nski set (denoted S)
is in exactly the classes indicated in figure \ref{cshl} of
the introduction.

\subsec{The generic Lusin set H}

The fact that no Lusin set
satisfies $\ufin(\ga,\ga)$ follows from Theorem \ref{gam6}.

\begin{theorem}[CH]\label{lusin} There exists a
Lusin set $H$ which is $\sone(\om,\om)$.
\end{theorem}
\pf
To construct a $\sone(\om,\om)$ Lusin set in the reals
enumerate all countable sequences of countable open
families as
$\{ ({\cal U}^{\beta  }_{n})_{n<\omega }:\beta  <\omega_{1}\} $.
Also enumerate all dense open subsets of the reals as
$(D_{\alpha })_{\alpha <\omega_{1}}$. We construct
$X$ recursively
as $\{ x_{\beta  }:\beta  <\omega _{1}\} $ as follows.
At stage
$\alpha$ of the construction we have
$$\{ x_{\beta  }:\beta <\alpha  \}
\mbox{ and }\{ (U_{n}^{\beta })_{n<\omega }:\beta <\alpha \}$$
satisfying for each $\beta <\alpha $
\begin{itemize}
\item[(i)] $x_\beta\in\cap \{D_\delta:\delta<\beta\}$,
\item[(ii)] $\{U_{n}^{\beta}: n<\omega \}$
is an $\omega$-cover of $\{ x_\delta:\delta <\alpha \}$,
\item[(iii)] if $({\cal U}^{\beta}_{n})_{n<\omega }$ was an
$\omega$-cover of $\{x_\delta:\delta<\beta\}$, then
$U_{n}^{\beta}\in {\cal U}^{\beta}_{n}$ for every $n$.
\end{itemize}

To see how to choose $x_{\alpha }$ and $(U_{n}^{\alpha })_{n<\omega}$
consider the $\alpha $'th sequence of open families: if
$({\cal U}^{\alpha }_{n})_{n<\omega }$ is a sequence of $\omega$-covers of
$\{x_{\beta} :\beta <\alpha  \}$ first extract an $\omega$-cover
$(U_{n}^{\alpha })_{n<\omega}$ so that
$U_{n}^{\alpha }\in {\cal U}^{\alpha }_{n}$ for each
$n<\omega $ (countable sets are
$\sone(\om,\om)$).
If $({\cal U}^{\alpha }_{n})_{n<\omega }$ is not a sequence of
$\omega$-covers
 of $\{ x_{\beta  }:\beta <\alpha  \}$ let
$U_{n}^{\alpha }={\Bbb R}$ for each $n<\omega$.  (${\Bbb R}$ is the set
of real numbers.)

Enumerate the finite subsets of $\{ x_{\beta  }:\beta <\alpha  \}$
as $\{ A_{k}:k<\omega \} $. For each $k$ and each $\beta \leq \alpha $
let
$$O_{k,\beta }=\bigcup\{ U_{n}^{\beta }:A_{k}\subseteq
U_{n}^{\beta }\} .$$
Then $O_{k,\beta }$ is dense and open. We
choose
$$
x_{\alpha }\in \bigcap_{\beta \leq \alpha  }D_{\beta }\cap
\bigcap_{k<\omega ,\beta \leq \alpha } O_{k,\beta }
$$
different from all $x_\beta$ with $\beta<\alpha$.
To see that $(U_{n}^{\beta })_{n<\omega }$ is an $\omega$-cover of
$\{ x_{\beta }:\beta \leq \alpha \}$ for each $\beta \leq \alpha$
it suffices to show that each $A_{k}\cup\{ x_{\alpha }\} $ is covered
by some $U_{n}^{\beta }$ for some $n<\omega$.  But
$x_{\alpha }\in O_{k,\beta }$ implies that there is an $n$ such that
$x_{\alpha}\in U_{n}^{\beta }$ and $A_{k}\subseteq U_{n}^{\beta }$.
We let $H=\{x_{\beta}: \beta<\omega_1 \}$.  To see that $H$
is $\sone(\om,\om)$, fix a sequence of $\omega$-covers
$({\cal U}_n)_{n<\omega}$. There is an $\alpha$ such that
$({\cal U}_n)_{n<\omega}=({\cal U}_n^{\alpha})_{n<\omega}$. Then at stage
$\alpha$ of the construction we extracted an appropriate
$\omega$-cover of $\{x_{\beta}:\beta\leq\alpha\}$ and inductive
hypothesis (ii) assures that it is also an $\omega$-cover of $H$.

\qed

The proof of Theorem \ref{lusin} only requires
that the covering number of the meager
ideal is equal to the continuum ($\covmeag={\goth c}$).
This requirement is equivalent to $MA$ for countable posets.
Adding Cohen reals over any model yields an $\sone(\om,\om)$ Lusin set
and hence our name generic Lusin set.

\section{Preservation of the properties.}\label{pres}

 Each of the properties in the diagram is inherited by closed
subsets and continuous images. The preservation theory is more
complicated for other topological constructions.

\begin{theorem}\label{contimage}  Let ${\sf G}$ be one of
 $\sone$, $\sfin$, or $\ufin$ and let ${\cal A}$ and
 ${\cal B}$ range over the set $\{\op,\om,\lm,\ga\}$.
 If $X$ has property ${\sf G}({\cal A},{\cal B})$ and
 $C$ is a closed subset of $X$, then $C$ has property
 ${\sf G}({\cal A},{\cal B})$.   If $f:X\to Y$ is continuous and onto and
 $X$ has the property ${\sf G}({\cal A},{\cal B})$, then so does $Y$.
\end{theorem}
\pf
The closure under taking closed subspaces is clear since if
$\cal U$ is a cover of $C$ in one of the classes $\{\op,\om,\lm,\ga\}$
for $C$, then
$${\cal V}=\{U\cup (X\setminus C): U\in {\cal U}\}$$
is in the same class for $X$.

To prove the closure under continuous image use that if
$\cal U$ is a cover of $Y$ in one of the classes $\{\op,\om,\lm,\ga\}$
for $Y$, then
$${\cal V}=\{f^{-1}(U): U\in {\cal U}\}$$
is in the same class for $X$.
\qed

\subsec{Finite powers}

We show that the classes $\sone(\om,\om)$, $\sfin(\om,\om)$,
and $\sone(\om,\ga)$ are the only ones closed under finite
powers.

\begin{lemma} \label{pow1} Let $X$ be a space and let
$n$ be a positive integer.
If ${\cal U}$ is an $\omega$--cover of $X$, then
$\{U^n:U\in{\cal U}\}$ is an $\omega$--cover of $X^n$.
\end{lemma}

\pf Observe that if $F$ is a finite subset of $X^n$, then there
is a finite subset $G$ of $X$ such that $F\subset G^n$.
\qed

\begin{lemma}\label{omcovpowers} Let $X$ be a topological space and
let $n$ be a positive integer. If ${\cal U}$ is an
$\omega$--cover for $X^n$, then there is an $\omega$--cover
${\cal V}$ of $X$ such that the open cover
$\{V^n:V\in{\cal V}\}$ of $X^n$ refines ${\cal U}$.
\end{lemma}

\pf Let ${\cal U}$ be an $\omega$--cover of $X^n$. Let $F$ be a
finite subset of $X$. Then $F^n$ is a finite subset
   of $X^n$. Since ${\cal U}$ is an $\omega$--cover of $X$, choose
an open set $U\in {\cal U}$ such that $F^n\subset U$.
   For any $n$--tuple $(x_1,\dots,x_n)$ in $F^n$, find for each
$i\in\{1,\dots,n\}$ an open set $U_i(x_1,\dots,x_n)\subset X$
    such that $x_i\in U_i(x_1,\dots,x_n)$, and
   $\prod_{i=1}^nU_i(x_1,\dots,x_n)\subset U$. Then, for each $x$
in $F$, let $U_x$ be the intersection of all the
   $U_i(x_1,\dots,x_n)$ which have $x$ as an element. Finally,
choose $V_F$ to be the set $\cup_{x\in F}U_x$, an open
   subset of $X$ which contains $F$, and which has the property
that $F^n\subset V^n_F\subset U$. Put
$${\cal V}= \{V_F:F \in[X]^{<\omega}\}.$$
Then ${\cal V}$ is as required.
\qed

While Lemma \ref{pow1} is true of $\gamma$--covers,
Lemma~\ref{omcovpowers} is not.

\begin{theorem}\label{s1omegprod} Let $n$ be a positive integer. If
a space $X$ has property $\sone(\om,\om)$,
so does $X^n$.
\end{theorem}
\pf Let $n$ be a positive integer and let $({\cal
U}_m:m=1,2,3,\dots)$ be a sequence of $\omega$--covers of $X^n$. 
By Lemma \ref{omcovpowers} for each $m$, we can choose
${\cal V}_m$ an $\omega$--cover of $X$ such
that $$\{V^n:V\in {\cal V}_m\}$$ is an
$\omega$--cover of $X^n$ which refines ${\cal U}_m$.

   Now apply the fact that $X$ is in $\sone(\om,\om)$ to
select from each ${\cal V}_m$ a set $V_m$ such
   that $\{V_m:m=1,2,3,\dots\}$ is an $\omega$--cover of $X$. Then,
since for each $m$ the set $\{V^n:V\in{\cal V}_m\}$
refines ${\cal U}_m$, we see that we can select from each
${\cal U}_m$ a set $U_m$ such that $V^n_m\subseteq U_m$.
But then the set $\{U_n:n=1,2,3,\dots\}$ is an $\omega$--cover
for $X$.
\qed

\begin{theorem}\label{sfinomegprod} Let $n$ be a positive integer
and let $X$ be a space. If $X$ has property
$\sfin(\om,\om)$, then $X^n$ also has this property.
\end{theorem}
\pf Let $({\cal U}_m:m=1,2,3,\dots)$ be a sequence of
$\omega$--covers of $X^n$. For each $m$, choose an $\omega$--cover
   ${\cal V}_m$ of $X$ such that $\{V^n:V\in{\cal V}_m\}$ refines
${\cal U}_m$. Now apply the fact that $X$ satisfies
   $\sfin(\om,\om)$: For each $m$ we find a finite subset
${\cal W}_m$ of ${\cal V}_m$ such that the collection
   $\cup_{m=1}^{\infty}{\cal W}_m$ is an $\omega$--cover of $X$.
For each $m$, choose a finite subset ${\cal Z}_m$ of
   ${\cal U}_m$ such that there is for each $W$ in ${\cal W}_m$ a
$Z$ in ${\cal Z}_m$ such that $W^n\subseteq Z$.
   Then $\cup_{m=1}^{\infty}{\cal Z}_m$ is an $\omega$--cover of
$X^n$.
\qed

\begin{theorem} Let $n$ be a positive integer
and let $X$ be a space. If $X$ has property
$\sfin(\om,\ga)$, then $X^n$ also has this property.
\end{theorem}
\pf
This is similar to the last two proofs.
\qed

\begin{theorem}
  [CH] None of the other classes (see figure \ref{cshl}) are
  closed under finite powers.
\end{theorem}
\pf
Note the examples L and S are such that
there sum $L+L$ and $S+S$ are homeomorphic to
the irrationals.

The function $\phi$ from $L\times L$ which assigns
to $(x,y)$ the point $\phi(x,y) = x+y$ is continuous.
But the space of irrationals does not
have property $\ufin(\ga,\op)$.
Since $\ufin(\ga,\op)$ is closed under continuous images
(see Theorem \ref{contimage})
$L\times L$ does not have property $\ufin(\ga,\op)$.

Similarly, $S \times S$ does not have property $\ufin(\ga,\op)$.
So none of the classes
containing either one of them is closed under finite powers.

\qed

We have seen that the inclusion
$\sone(\om,\om)\subseteq\sone(\op,\op)$ may be
proper, e.g. the special Lusin set L is in $\sone(\op,\op)$
but not in $\sone(\om,\om)$.
We now give an important fact about these two classes, which
characterizes
$\sone(\om,\om)$ as a subset of $\sone(\op,\op)$.

\begin{theorem} Let $X$ be a space. Then the following are
equivalent: \label{sonepow}
\begin{enumerate}
\item $X$ satisfies $\sone(\om,\om)$.
\item Every finite power of $X$ satisfies $\sone(\op,\op)$
      (Rothberger property $C^{\prime\prime}$).
\end{enumerate}
\end{theorem}
\pf The implication $1\Rightarrow 2$ follows immediately
from Theorem \ref{s1omegprod} and the fact that
$\sone(\om,\om)$ is a subclass of $\sone(\op,\op)$.

\medskip

The implication $2\Rightarrow 1$ is proven as follows: Let
$({\cal U}_n:n\in\omega)$ be a sequence of $\omega$--covers
   of $X$. Write the set of positive integers as a union of
countably many disjoint infinite sets, say
   $Y_1, Y_2,\dots, Y_n,\dots$. For each $m$ and for each $k$ in
$Y_m$ put ${\cal V}_k = \{U^m:U\in{\cal U}_k\}$.

   Then by Lemma \ref{pow1}, for each $m$ the sequence
$({\cal V}_k:k\in Y_m)$ is a
sequence of $\omega$--covers of $X^m$. By Hypothesis $2$
   we find for each $m$ a sequence
$$(U^m_k:k\in Y_m)$$
such that for each $k$, $U_k\in{\cal U}_k$, and such that
   $\{U^m_k:k\in\omega\}$ is an open cover of $X^m$.

   The sequence $(U_k:k\in\omega)$ is an
$\omega$--cover of $X$. For let $F$ be a finite subset
   of $X$, say $F=\{x_1,\dots,x_m\}$, enumerated bijectively. Then
$(x_1,\dots,x_m)$ is an element of $X^m$. Thus, choose a
   $k$ in $Y_m$ such that $(x_1,\dots,x_m)$ is in $U^m_k$. Then $F$
is a subset of $U_k$.
\qed

The Borel Conjecture, that every strong measure zero set is
countable, implies that the two classes $\sone(\om,\om)$ and $\sone(\op,\op)$
coincide. The Borel Conjecture was proved consistent by Laver.

\begin{problem} Is it true that if there is an uncountable set of
 real numbers which has property $\sone(\om,\om)$, then there is a
 set of real numbers which has property $\sone(\op,\op)$ but does not have
 property $\sone(\om,\om)$?
\end{problem}

\medskip

We shall now prove the analogue of Theorem \ref{sonepow} for
$\sfin(\op,\op)$ and $\sfin(\om,\om)$.

\begin{theorem} \label{finitepower} For a space $X$ the following
are equivalent:
\begin{enumerate}
 \item Every finite power of $X$ has property
       $\sfin(\op,\op)$.
 \item $X$ has property $\sfin(\om,\om)$.
\end{enumerate}
\end{theorem}
\pf The implication $2\Rightarrow 1$: This follows from
Theorem \ref{sfinomegprod}.

\medskip

   We now work on the implication $1\Rightarrow 2$: Let $({\cal
U}_n:n\in\omega)$ be a sequence of
   $\omega$--covers of $X$. Let $(Y_k:k\in\omega)$ be a pairwise
disjoint sequence of infinite sets
   of positive integers whose union is the set of positive
integers. For each $m$, for each $k$ in 
   $Y_m$, put ${\cal V}_k = \{U^m:U\in{\cal U}_k\}$.
   Then for each $m$ by Lemma \ref{omcovpowers}, the sequence
$({\cal V}_k:k\in Y_m)$ is a sequence of $\omega$--covers of $X^m$.

   Applying $1$ for each $m$, we find for each $m$ a sequence
$({\cal W}_k:k\in Y_m)$ such that
\begin{itemize}
\item{for each $k\in Y_m$, ${\cal W}_k$ is a finite subset of
${\cal U}_k$, and}
\item{$\cup_{k\in Y_m}\{U^m:U\in{\cal W}_k\}$ is an open cover of
$X^m$.}
\end{itemize}
But then $\cup_{k=1}^{\infty}{\cal W}_k$ is an $\omega$--cover
of $X$.
\qed

None of our classes are closed under finite products.
Todorcevic \cite{To} showed that there exist
two (nonmetrizable) topological spaces $X$ and $Y$ that satisfy
$\sone(\om,\ga)$ ($\gamma$-set), but whose product
does not satisfy $\ufin(\ga,\op)$ (Menger).
Thus none of our properties are closed under finite products.

If we restrict our attention to separable metric spaces it
also is the case assuming CH that none of our classes are closed
under finite products.  For the class $\sone(\om,\ga)$ note
that Galvin-Miller \cite{G-M} using a result of Todorcevic
showed that there are $\gamma$-sets whose product is not
a $\gamma$-set.  For the classes $\sone(\om,\om)$ and
$\sfin(\om,\om)$ construct a pair of generic Lusin sets
$H_0$ and $H_1$ such that $H_0+H_1=\zz$

\medskip

\noindent {\bf Remark}
The special Lusin set L gives a partial
answer to a problem of Lelek (see \cite{Le}).
It shows that it is relatively consistent with ZFC
that there exists a separable metrizable space $L$ that has
property $\ufin(\ga,\op)$, but
does not have property $\ufin(\ga,\op)$ in each finite power.
In Lelek, $\ufin(\ga,\op)$ is referred to as the ``Hurewicz
property'' in contrast to our notation as the Menger property.

\medskip

\noindent {\bf Remark} It is relatively consistent with ZFC that
for every $n\geq 1$ there exists a separable metric space $X$ such
that $X^n$ has property $\ufin(\ga,\op)$ but $X^{n+1}$ does not
have property $\ufin(\ga,\op)$
(see Just \cite{J1} and Stamp \cite{WS}).

\medskip

\noindent {\bf Remark}
It was shown in Just \cite{WJ} that preservation of $\ufin(\ga,\om)$
under direct sums is independent of ZFC.

\subsec{Finite or countable unions}

It is well-known and easy to prove that each of the classes
\begin{itemize}
\item $\sfin(\op,\op)$ (Rothberger property $C^{\prime\prime}$),
\item $\ufin(\ga,\ga)$ (Hurewicz property), and
\item $\ufin(\ga,\op)$ (Menger property)
\end{itemize}
are closed under taking countable unions.
It also easy to prove that $\sone(\ga,\lm)$ is closed under taking countable
unions. The class $\sone(\om,\ga)$ (Gerlits-Nagy property $\gamma$-sets)
is not closed under taking  finite unions
(see Galvin-Miller \cite{G-N}).

\begin{problem}
  Which of the remaining classes are closed under taking
  finite or countable unions?
\end{problem}

\section{Cardinal equivalents}\label{card}

\begin{figure}
\unitlength=.95mm
\begin{picture}(140.00,100.00)(10,10)
\put(20.00,20.00){\makebox(0,0)[cc]
{\shortstack {$\sone(\om,\ga)$\\ ${\goth p}$ } }}
\put(60.00,20.00){\makebox(0,0)[cc]
{\shortstack {$\sone(\om,\om)$\\ ${\covmeag}$ } }}
\put(100.00,20.00){\makebox(0,0)[cc]
{\shortstack {$\sone(\op,\op)$\\ ${\covmeag}$ } }}
\put(20.00,60.00){\makebox(0,0)[cc]
{\shortstack {$\sone(\ga,\ga)$\\ ${\goth b}$ } }}
\put(60.00,60.00){\makebox(0,0)[cc]
{\shortstack {$\sone(\ga,\om)$\\ ${\goth d}$ } }}
\put(100.00,60.00){\makebox(0,0)[cc]
{\shortstack {$\sone(\ga,\lm)$\\ ${\goth d}$ } }}
\put(80.00,40.00){\makebox(0,0)[cc]
{\shortstack {$\sfin(\om,\om)$\\ ${\goth d}$ } }}
\put(80.00,80.00){\makebox(0,0)[cc]
{\shortstack {$\sfin(\ga,\om)$\\ ${\goth d}$ } }}
\put(60.00,100.00){\makebox(0,0)[cc]
{\shortstack {$\ufin(\ga,\ga)$\\ ${\goth b}$ } }}
\put(100.00,100.00){\makebox(0,0)[cc]
{\shortstack {$\ufin(\ga,\om)$\\ ${\goth d}$ } }}
\put(140.00,100.00){\makebox(0,0)[cc]
{\shortstack {$\ufin(\ga,\op)$\\ ${\goth d}$ } }}

\put(70.00,100.00){\vector(1,0){20.00}}
\put(110.00,100.00){\vector(1,0){20.00}}
\put(86.00,85.00){\vector(1,1){10.00}}
\put(25.00,65.00){\vector(1,1){30.00}}
\put(105.00,65.00){\vector(1,1){30.00}}
\put(65.00,65.00){\vector(1,1){10.00}}
\put(28.00,61.00){\vector(1,0){20.00}}
\put(70.00,61.00){\vector(1,0){20.00}}
\put(80.00,45.00){\line(0,1){12.00}}
\put(80.00,64.00){\vector(0,1){11.00}}
\put(20.00,25.00){\vector(0,1){29.00}}
\put(60.00,25.00){\vector(0,1){29.00}}
\put(64.00,27.00){\vector(1,1){10.00}}
\put(100.00,25.00){\vector(0,1){29.00}}
\put(28.00,20.00){\vector(1,0){20.00}}
\put(70.00,20.00){\vector(1,0){20.00}}
\end{picture}
\caption{Cardinals $\non(P)$ \label{cardfig}}
\end{figure}
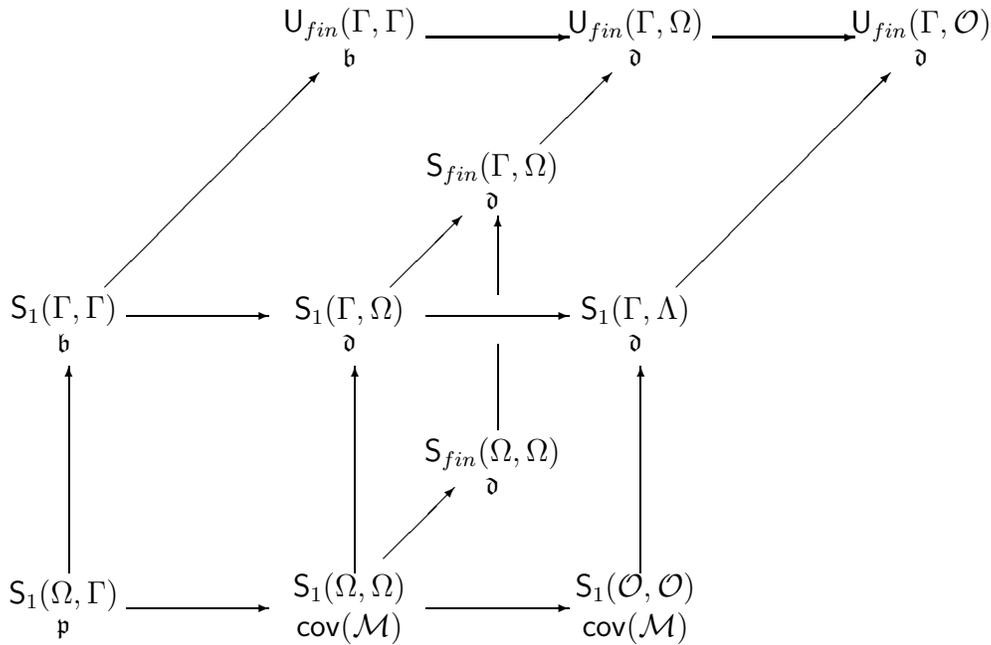

We now consider the connection between the properties and some
well known cardinal invariants of $P(\omega)/Fin$. See
Vaughan \cite{vaughan}  for the definitions, but briefly

\smallskip

\noindent ${\goth p}$ is least cardinality of a family of
sets in $[\omega]^\omega$
with the finite intersection property but no pseudo intersection,

\smallskip

\noindent {\goth d} is the minimal cardinality of a dominating family
in $\oo$,

\smallskip

\noindent ${\goth b}$ the minimal cardinality of an unbounded family, and

\smallskip

\noindent $\covmeag$ is the minimal cardinality of a covering of the
real line by meager sets.

\smallskip

In particular,
if $P$ is one of the eleven properties in the diagram
(figure \ref{cshl}) or is one of
the splitting properties $\split(\om,\om )$ or
$\split(\lm,\lm )$ we will determine:

\medskip

$non(P)$:=

the minimum cardinality of a set of
reals that fails to have property $P$.

\medskip
\noindent Note obviously that if $P\rightarrow Q$, then
$\non(P)\leq\non(Q)$.
Some of these cardinals are well known and we simply state the results
and refer the reader to the appropriate references.

\begin{theorem}\label{gam0} (Galvin-Miller \cite{G-M})
$\non(\sone(\om,\ga))={\goth p}$.
\end{theorem}

\begin{theorem}\label{gam2} (Fremlin-Miller \cite{F-M})
$\non(\sone(\op,\op))=\covmeag$.
\end{theorem}

The next two theorems are due to Hurewicz and they imply that the
minimum cardinality of a set of reals that is not $\ufin(\ga,\ga)$ is
${\goth b}$, and the minimum cardinality of a set of reals that is
not $\ufin(\ga,\op)$ is ${\goth d}$.

\begin{theorem}\label{gam6} (Hurewicz \cite{Hu2})
A set $X$ is $\ufin(\ga,\ga)$ if and
only if every continuous image of $X$ in $\oo$ is bounded.
Hence $\non(\ufin(\ga,\ga))={\goth b}$.
\end{theorem}

\begin{theorem}\label{gam8} (Hurewicz \cite{Hu2})
A set $X$ is $\ufin(\ga,\op)$ if and
only if every continuous image of $X$ in $\oo$ is not
dominating. Hence $non(\ufin(\ga,\op))={\goth d}$.
\end{theorem}

Next we determine $\non(P)$ for all the other properties
in figure \ref{cshl} in the introduction.

\begin{theorem}\label{gam4} $\non(\sone(\ga,\om))={\goth d}$.
\end{theorem}

\pf
Since $\sone(\ga,\om)\subseteq \ufin(\ga,\op)$ we have
that $$\non(\sone(\ga,\om))\leq \non(\ufin(\ga,\op)).$$
Also by Theorem \ref{gam8} we have $\non(\ufin(\ga,\op))={\goth d}$,
so $\non(\sone(\ga,\om))\leq {\goth d}$.

Conversely, suppose that $X$ is a set of reals that fails to be
$\sone(\ga,\om)$. Fix a sequence of $\gamma$-covers
$({\cal U}_{n})_{n\in \omega}$ witnessing the failure of
$S_{1}(\ga ,\om ).$ Fix an enumeration of each cover
${\cal U}_{n}=\{U_{n}^{i}:i\in \omega\}$. For each finite set
$F\subseteq X$ define $f_{F}\in \oo$ by
$$
f_{F}(n)=min\{i:\forall j > i,\; F\subseteq U_{n}^{i}\}.
$$
As each ${\cal U}_{n}$ is a $\gamma$-cover if $i>f_{F}(n)$, then
$F\subseteq U_{n}^{i}$. Therefore, $$\{f_{F}:F\in[X]^{<\omega}\}$$
must be a dominating family. Otherwise there is a $g$ not dominated
by any such $f_{F}$. I.e., for each finite $F\subseteq X$, there is
an integer $n$ such that $g(n)>f_{F}(n)$. This implies that
$\{U_{n}^{g(n)}:n\in \omega\}$ is an $\omega$-cover,
contradicting the failure of $S_{1}(\ga ,\om)$.
So $\non(\sone(\ga,\om))\geq {\goth d}$.
\qed

\begin{theorem}\label{sfinomega}
$\non(\sfin(\om,\om))={\goth d}$.
\end{theorem}

\pf Identical to the proof of \ref{gam4}.  One only needs to modify
the definition of $f_F$ to
$$f_{F}(n)=min\{i:F\subseteq U_{n}^{i}\}$$
and take ${\cal V}_n=\{U^i_n:i\leq g(n)\}$.
\qed

\begin{theorem}\label{gam3} $\non(\sone(\ga,\ga))={\goth b}$.
\end{theorem}
\pf

Using $\sone(\ga,\ga)\subseteq \ufin(\ga,\ga)$ and Theorem~\ref{gam6}
it follows that $$\non(\sone(\ga,\ga))\leq {\goth b}.$$

Conversely, suppose that $X$ is a set of reals and that
$({\cal U}_{n})_{n\in \omega}$ is sequence of $\gamma $-covers
witnessing the failure of $S_{1}(\ga ,\ga )$. For each
$x\in X$ define $f_{x}\in \oo$ by
$$f_{x}(n)=min\{i:\forall j\geq i, x\in U_{n}^{j}\}.$$
If $g$ were to dominate each $f_{x}$,
then $(U_{n}^{g(n)})_{n\in \omega}$ would be a $\gamma$-cover, a
contradiction. Therefore $\{f_{x}:x\in X\}$ is an unbounded
family. Hence $\non(\sone(\ga,\ga))\geq {\goth b}$.
\qed

\begin{theorem}\label{gam1}
$\non(\sone(\om,\om))=\covmeag$.
\end{theorem}

\pf
The inclusion $\sone(\om,\om)\subseteq \sone(\op,\op)$ and
Theorem~\ref{gam2} give us the inequality
$\non(\sone(\om,\om))\leq\covmeag$.

Conversely fix $X$ a set of reals and
$({\cal U}_{n})_{n\in \omega}$
a sequence of $\omega$-covers witnessing
the failure of $\sone(\om,\om)$.  For each finite
$F\subseteq X$ let
$$
K_{F}=\{f\in \oo :(\forall n\in \omega)(F\not\subset U_{n}^{f(n)}\}.
$$
Since for each $f\in \oo$ there is a finite
$F\subseteq X$ such that $F\not\subset U_{n}^{f(n)}$, we have that
$\oo=\bigcup\{K_{F}:F\in [X]^{<\omega}\}$.
Furthermore each $K_{F}$ is closed and
nowhere dense.  Hence $\non(\sone(\om,\om))\geq\covmeag$.
\qed

Our results are summarized in figure \ref{cardfig}.
Classical results about the relationships between the cardinals
${\goth p}$, ${\goth b}$, ${\goth d}$ and $cov({\cal M})$ give alternative
proofs that many of the implications in our diagram cannot be
reversed.

\subsec{$\split(\lm,\lm)$ and $\split(\om,\om)$}

These properties were defined in
Scheepers \cite{S}: for classes of covers ${\cal A}$ and
${\cal B}$, a space has
property $\split({\cal A},{\cal B})$ iff every open cover
${\cal U\in A}$ can be partitioned into two subcovers ${\cal U}_{0}$ and
${\cal U}_{1}$ both in ${\cal B}$.
Recall that a family ${\cal R}\subseteq [\omega]^{\omega}$
is said to be a {\em reaping family} if for
each $x\in [\omega]^{\omega}$ there is a $y\in {\cal R}$ such
that either $y\subseteq^{*}x$ or $y\subseteq^{*} \omega\setminus x$.
The minimal cardinality of a reaping family is denoted by
${\goth r}$, and the minimal cardinality of a base for a nonprincipal
ultrafilter is denoted by ${\goth u}$.

\begin{theorem}\label{large1}
$\non(\split(\lm,\lm))={\goth r}$.
\end{theorem}

\pf Suppose that $X\subseteq[\omega]^{\omega}$ is a
reaping family. Consider the open family
$${\cal U}=\{B^{1}_{n}:n\in \omega\}$$
where
$$B^{1}_{n}=\{x\in [\omega]^{\omega}:n\in x\}
\mbox{ and } B^{0}_{n}=\{x\in [\omega]^{\omega}:n\not\in x\}.$$
Clearly ${\cal U}$ is a large cover of any subset
of $[\omega]^{\omega}$.  We will often refer to it as the
canonical large cover. Since $X$ is a reaping family, this cover
cannot be partitioned into two large subcovers.
Conversely, suppose that $X$ is a set of reals and
${\cal U}=\{U_{n}:n\in \omega\}$ is a large cover of $X$.
For each $x\in X$ let
$$
A_{x}=\{n\in \omega:x\in U_{n}\}.
$$
If ${\cal F}$ is the collection of all such $A_{x}$'s, then
${\cal F}$ is a reaping family. For if $A\subseteq \omega$ is such that
for all $x\in X$ both $A_{x}\cap A$ and $A_{x}\setminus A$ are
infinite, then
$$\{U_{n}:n\in A\}\cup\{U_{n}:n\not\in A\}$$
is a splitting of ${\cal U}$ into disjoint large subcovers.
\qed

The proof yields a bit more.

\begin{theorem}\label{large2} A set of reals $X$ is $\split (\lm ,\lm )$
with respect to clopen covers if and only if every
continuous image of $X$ in $[\omega]^{\omega}$ is not a reaping
family.
\end{theorem}

\pf Suppose that $X$ is a set of reals,
$f:X\rightarrow [\omega]^{\omega}$ is continuous and that $f(X)$
is a reaping family. The canonical large cover is in fact a clopen
family. Therefore the collection $f^{-1}({\cal
U})=\{f^{-1}(B^{1}_{n}):n\in \omega\}$ is a large clopen cover of
$X$. Suppose $f^{-1}({\cal U})={\cal V}_{0}\cup{\cal V}_{1}$ is a
partition. Then we have the corresponding partition of
$\omega=A_{0}\cup A_{1}$ where
${\cal V}_{i}=\{f^{-1}(U_{n}):n\in A_{i}\}$.
As $f(X)$ is a reaping family, there is an $x\in X$ such
that for either $i=0$ or $1$, $f(x)\subseteq^{*} A_{i}$. Then
${\cal V}_{i}$ is not large at $x$.
Therefore $X$ is not $\split (\lm ,\lm )$ with respect to
the clopen cover $f^{-1}({\cal U})$.

\smallskip
\noindent Conversely, suppose that $X$ is not
$\split(\lm,\lm ) $ with respect to some large clopen cover
${\cal U}=\{U_{n}:n\in \omega\}$.  For each $x\in X$ define
$f_{x}\in [\omega]^{\omega}$
by $n\in f_{x}$ iff $x\in U_{n}$. Since $U$ is large, each $f_{x}$
is infinite. As above, since ${\cal U}$ cannot be split,
$\{f_{x}:x\in X\}$ is a reaping family. Therefore it suffices to
check that the mapping $f:x\rightarrow f_{x}$ is continuous. But
the collection of $\{B_n^i:n\in \omega,i=0,1\}$ forms a subbase
for $[\omega]^{\omega}$,
and clearly $f^{-1}(B_{n}^{1})=U_{n}$ and
$f^{-1}(B_{n}^{0})=X\setminus
U_{n}$ therefore $f$ is continuous (this is the only place where we
need the restriction to clopen covers).
\qed

\begin{theorem}\label{omega}  $\non(\split(\om,\om))={\goth u}$.
\end{theorem}

\pf Suppose that $X\subseteq[\omega]^{\omega}$ is a filter-base.
Then the canonical large cover in $[\omega]^{\omega}$
is in fact an $\omega$-cover of $X$. If $X$ is a base for an
ultrafilter, then the
canonical cover cannot be partitioned into two $\omega$-subcovers.

\smallskip
\noindent Conversely, suppose that $X$ is a set of reals and ${\cal U}$
is an $\omega$-cover of $X$. For each $x\in X$ let 
$$
{\cal U}_{x}=\{U\in {\cal U}:x\in U\}.
$$
If ${\cal F}$ is the collection of all such ${\cal U}_{x}$'s, then
${\cal F}$ forms a filterbase
on ${\cal U}$ and if ${\cal U}$ cannot be split into two
$\omega$-covers, then ${\cal F}$ generates a nonprincipal ultrafilter.
\qed

Analogously to Theorem \ref{large2} we can prove:

\begin{theorem}\label{omega2} A set of reals $X$ is
$\split(\om,\om)$ with respect to clopen covers if and only if every
continuous image of $X$ in $[\omega]^{\omega}$ does not
generate an ultrafilter.
\end{theorem}

Note that a base for an ultrafilter is a reaping family, and
therefore ${\goth r}\leq{\goth u}$. In Bell-Kunen \cite {belku}
it is proven consistent that
this inequality may be strict. Therefore
$\split(\lm,\lm )\not\Rightarrow \split (\om ,\om )$. Similarly
neither ${\goth r}$ nor ${\goth u}$ are comparable to ${\goth d}$,
therefore there are no implications between either
$\split(\lm,\lm )$ or $\split(\om,\om)$ and any of
the six classes in figure \ref{cardfig} whose `$\non$' is
equivalent to ${\goth d}$.
In Scheepers \cite{S} it is shown that
\begin{itemize}
  \item $\ufin(\ga,\ga)\Rightarrow \split(\lm,\lm)$ (Cor 29), and
  \item $\sone(\op,\op)\Rightarrow \split(\lm,\lm)$ (Thm 15).
\end{itemize}
Note that while both ${\goth b}\leq {\goth r}$
and $\covmeag\leq {\goth r}$, it is consistent that these
inequalities are strict (see Vaughan \cite{vaughan}).
So neither of these implications can be reversed.
\begin{problem}
  Does $\split(\om,\om)\Rightarrow \split(\lm,\lm)$?
\end{problem}

\section{The Hurewicz Conjecture and the Borel Conjecture.}\label{hur}

   Every $\sigma$--compact space belongs to $\ufin(\ga,\ga)$. It is
also well-known that not every space belonging to $\ufin(\ga,\ga)$ need
be $\sigma$--compact.  We now look at the traditional examples of sets
of reals belonging
to $\ufin(\ga,\ga)$, and show that some of these belong to
$\sone(\ga,\ga)$, while others do not. Since $\sone(\ga,\ga)$ is
contained in $\sone(\ga,\lm)$,
and the unit interval is not an element of $\sone(\ga,\lm)$, we see that
the $\sigma$--compact spaces do not in general belong to the class
$\sone(\ga,\ga)$.

   On page 200 of \cite{Hu}, W. Hurewicz conjectures:

\medskip

\noindent
{\em
[Hurewicz] A set of real numbers has property
$\ufin(\ga,\ga)$ if, and only if, it is
$\sigma$--compact.\footnote{``Es
entsteht nun  die Vermutung dass durch die (warscheinlich
sch\"arfere)
Eigenschaft $E^{**}$ die halbkompakten Mengen $F_{\sigma}$
allgemein charakterisiert sind.''}
}

\medskip

   The existence of a Sierpi\'nski set violates this conjecture. As
we have seen earlier, Sierpi\'nski sets are elements of
$\sone(\ga,\ga)$.

The following result shows that Hurewicz's conjecture fails
in ZFC.
\begin{theorem}\label{nothc}
  There exists a separable metric space $X$ such that
  $|X|=\omega_1$, $X$ is not $\sigma$-compact and $X$ has property
  $\ufin(\ga,\ga)$.
  This $X$ also has property $\sone(\ga,\om)$.
\end{theorem}
\pf

\medskip\noindent Case 1. ${\goth b}>\omega_1$.

In this case every $X$ of size $\omega_1$ is in $\sone(\ga,\ga)$,
hence in both $\ufin(\ga,\ga)$ and $\sone(\ga,\lm)$.
(In this case also in $\sfin(\Omega,\Omega)$.)

\medskip\noindent Case 2. ${\goth b}=\omega_1$.

In this case we will use a construction similar to one in
\cite{G-M}.
Build an $\omega_1$-sequence $\la x_\alpha:\alpha<\omega_1 \ra$ of
elements of $[\omega]^\omega$ such that $\alpha < \beta$
implies $x_\beta\subseteq^* x_\alpha$ and if $f_\alpha:\omega\to
x_\alpha$
is the increasing enumeration of $x_\alpha$, then
for every $g\in\oo$ there exists $\alpha$ such that
for infinitely many $n$ we have $g(n)<f_\alpha(n)$.

\begin{claim}\label{claim1} For any $S\in[\omega]^\omega$
there exists
$\alpha<\omega$ such that there exists infinitely many
$n$ such that $|[f_\alpha(n),f_\alpha(n+1))\cap S|\geq 2$.
\end{claim}
To prove Claim \ref{claim1} suppose not and let $g$
eventually
dominate all the increasing enumerations of sets $S^*$ such
that $S^*=^*S$.  Then $g$ eventually
dominates the $f_\alpha$'s, contradiction. This completes the proof of
Claim \ref{claim1}.

\begin{claim}\label{claim2} Let
$X=[\omega]^{<\omega}\cup\{x_\alpha:\alpha<\omega_1\}$.
Then for every sequence $\la{\cal U}_n:n\in\omega\ra$ of $\omega$--covers of
$X$ (or even just of $[\omega]^{<\omega}$) there exists an
$A\in [\omega]^\omega$, $\la V_n\in {\cal U}_n: n\in A\ra$ and
$\alpha<\omega_1$ such that for all $\beta\geq\alpha$ we
have $x_\beta\in V_n$ for all but finitely many $n\in A$.
\end{claim}
To prove Claim \ref{claim2} construct $k_n$, an increasing sequence
in
$\omega$, such that there exists $V_n\in{\cal U}_n$ with the
property that
$$\{\;x\subseteq \omega\; :\; x\cap (k_n,k_{n+1})=
\emptyset\;\}\subseteq V_n.$$
(Use that ${\cal U}_n$ is a $\omega$--cover to pick
$V_n\supseteq [k_n+1]^{<\omega}$ and then choose $k_{n+1}$.)
It follows from Claim \ref{claim1} that there exists an
$\alpha<\omega_1$,
$A\in[\omega]^\omega$, and an increasing sequence $\la m_n: n \in
A\ra$ such that
for every $n\in A$
$$\{\;x\subseteq \omega\;:\; x\cap
(f_{\alpha}(m_n),f_{\alpha}(m_n+1))
=\emptyset\;\}\subseteq V_n.$$
It follows  $x_\beta\in V_n$  for all $\beta\geq\alpha$ for all but
finitely many $n\in A$.  This completes the proof
of Claim \ref{claim2}

\medskip

Now we show that our set $X$ in this case is in both
$\ufin(\ga,\ga)$ and
$\sone(\om,\om)$ (and hence $\sone(\ga,\lm)$).  First we show that
it satisfies a property we might call $\sone(\ga,\ga)^*$.
\begin{quote}
  Given any $\la{\cal U}_n:n \in \omega\ra$ a sequence of
  $\gamma$--covers of $X$, there exists
  $\la V_n \in {\cal U}_n : n \in \omega\ra$ and a countable
$Y\subseteq X$
  such that $\la V_n: n \in \omega\ra$ is a $\gamma$--cover of
  $X\setminus Y$.
\end{quote}

If $\sfin(\ga,\ga)^*$ is defined analogously, then it is
easy to see using the same proof as Theorem \ref{sone_eq_sfin}
that $\sfin(\ga,\ga)^*$ is equivalent to   $\sone(\ga,\ga)^*$.
Clearly Claim \ref{claim2} implies $\sfin(\ga,\ga)^*$.

$\sone(\ga,\ga)^*$ implies $\ufin(\ga,\ga)$ because
we may first pick $\la V_n\in {\cal U}_n: n \in \omega\ra$ a
$\gamma$--cover of $X\setminus Y$ and then pick
$\la W_n\in {\cal U}_n: n \in \omega\ra$ a $\gamma$--cover of
$Y$.  Then $\la V_n \cup W_n: n \in \omega\ra$ is a
$\gamma$--cover of $X$.
\qed

To see that $X$ is in $\sone(\om,\om)$  we need the following
claim:

\begin{claim}\label{claim3} For every $B\in [\omega]^\omega$, sequence
$\la{\cal U}_n:n\in B\ra$ of $\omega$--covers of
$X$, and countable $Y\subseteq X$ there exist
$A\in [B]^\omega$, $\la V_n\in {\cal U}_n: n\in A\ra$ and
$Z\subseteq X$ countable
such that $Y$ and $Z$ are disjoint and
$\la V_n\in {\cal U}_n: n\in A\ra$ is a
$\gamma$ cover of $X\setminus Z$.
\end{claim}

\pf Let $Y=\{y_n:n\in\omega\}$ and
apply Claim \ref{claim2} to the $\omega$--covers
defined by
$${\cal U}^\prime_n =\{U\in {\cal U}_n: \{y_i :i< n\}\subseteq
U\}$$
for $n\in B$. This completes the proof of Claim \ref{claim3}.

Using Claim \ref{claim3}
for every sequence $\la{\cal U}_n:n\in\omega\ra$ of
$\omega$--covers of
$X$ inductively construct $A_i\in [\omega]^\omega$,
$\la V_n\in {\cal U}_n: n\in A_i\ra$ and $Y_i\subseteq X$ countable
such that
\begin{itemize}
\item $A_i\cap A_j=\emptyset$ for $i\not = j$.
\item $Y_i\cap Y_j=\emptyset$ for $i\not = j$.
\item $\la V_n\in {\cal U}_n: n\in A_i\ra$ is a $\gamma$--cover
of $X\setminus Y_i$.
\end{itemize}
\medskip
(At stage $n$ take $Y=\cup\{Y_i:i<n\}$ and
$B=\omega \setminus (\cup\{A_i:i<n\})$.  Apply Claim 3
and let $Y_n=Z$ and
cut down $A_n$, if necessary, to ensure that
$\cup\{A_i:i \leq n\}$ is
coinfinite.)

Since the $(Y_i:i<\omega)$ and the $(A_i:i<\omega)$ are pairwise
disjoint families, letting $A=\bigcup _{i\in\omega}A_i$, $(V_n:n\in
A)$
is an $\omega$-cover of $X$.   Hence $X$ has property
$\sone(\om,\om)$. This completes the proof of Theorem \ref{nothc}.
\qed

\begin{problem}
  Is the set $X$ constructed in Case 2 of Theorem \ref{nothc}
  a $\gamma$-set, i.e., $\sone(\om,\ga)$? 
\end{problem}

The Borel conjecture implies that every set
in $\sone(\op,\op)$ is countable (hence
every set in $\sone(\om,\om)$ or $\sone(\om,\ga)$ is
countable).
  Theorem \ref{nothc}
and the Cantor set along with the last example
rules out an analogous conjecture for all except
$\sone(\ga,\ga)$.  So we ask:

\begin{problem}
  Is it consistent, relative to the consistency of
  ${\sf ZF}$, that every set in $\sone(\ga,\ga)$ is
  countable?
\end{problem}

One may also ask if all the pathological examples of
sets having property $\ufin(\ga,\ga)$ occur because
of the presence of such sets in $\sone(\ga,\ga)$; here is one
formalization of this question.

\begin{problem} Let $X$ be a set of real numbers which does not
contain a perfect set of real numbers
but which does have the Hurewicz property. Does $X$ then belong
to $\sone(\ga,\ga)$?
\end{problem}

\subsec{$\ufin(\ga,\ga)$ and perfectly meager sets.}

We now prove a theorem which implies that the $\sone(\ga,\ga)$--sets 
are contained in another class of sets that were introduced in the early
parts of this century. Recall that a set $X$ of real numbers is
{\em perfectly meager} (also called ``always of first category") if,
for every perfect
set $P$ of real numbers, $X\cap P$ is meager in the relative topology of 
$P$.

 \begin{theorem} If a set of reals $X$ is in $\ufin(\ga,\ga)$ and
 contains no perfect subset, then $X$ is perfectly meager.
 \end{theorem}
 \pf
 Let $P$ be a perfect set of real numbers. Since $X$ contains no
 perfect set, $P\setminus X$ is a dense subset of $P$. Let $D$ be a countable
 dense subset of $P$ which is contained in $P\setminus X$, and enumerate $D$
 bijectively as $(d_n:n=1,2,3,...)$.

 Fix $k$. For each $x$ in $X$ choose open intervals $I^k_x$ and $J^k_x$
 such that
 \begin{enumerate}
 \item{$I^k_x$ is centered at $x$,}
 \item{$J^k_x$ is centered at $d_k$, and}
 \item{the closures of these intervals are disjoint.}
 \end{enumerate}
      Let $\{I^k_{x^k_n}:n=1,2,3,...\}$ be a countable subset of
      $\{I^k_x:x \in X\}$ which covers $X$. Then for each $n$ define
      $I^k_n = \cup_{j\leq n}I^k_{x^k_j}$, and $J^k_n = \cap_{j\leq
      n}J^k_{x^k_j}$. Then ${\cal U}_k=\{I^k_n:n=1,2,3,...\}$ is a
      $\gamma$--cover of $X$.

      Apply $\ufin(\ga,\ga)$ to the sequence $({\cal
      U}_k:k=1,2,3,..)$. For each $k$ we find an $n_k$ such that
      $(I^k_{n_k}:k=1,2,3,..)$ is a $\gamma$ cover for $X$.
      For each $j$ put $G_j = \cup_{k\geq j}J^k_{n_k}$. Then each $G_j\cap P$ 
      is a dense open subset of $P$ (as it contains all but a finite piece of 
      $D$). The intersection of these sets is a dense $G_{\delta}$ subset of
      $P$, and is disjoint from $X\cap P$.
      Thus, $X\cap P$ is a meager subset of $P$. \qed

 \begin{corollary} Every element of $\sone(\ga,\ga)$ is perfectly meager.
 \end{corollary}
 \pf  We have seen (Theorem \ref{thm17})
      that sets in $\sone(\ga,\ga)$ do not contain perfect
      sets of real numbers.
      But $\sone(\ga,\ga)\subseteq\ufin(\ga,\ga)$.
 \qed

 In Theorem 2 of Galvin-Miller \cite{G-M} it was shown that if
 a subset $X$ of the
 real line is in  $\sone(\om,\ga)$, then for every
 $G_{\delta}$ set $G$ which  contains
 $X$, there is an $F_{\sigma}$ set $F$ such that $X\subseteq F\subseteq
 G$.   In fact, this property characterizes $\ufin(\ga,\ga)$.

\begin{theorem}
  For a set $X$ of real numbers, the following are equivalent:
 \begin{enumerate}
 \item  $X$ has property $\ufin(\ga,\ga)$.
 \item  For every $G_{\delta}$--set $G$ which contains $X$, there is a
    $F_{\sigma}$--set $F$ such that $X\subseteq F\subseteq G$.
 \end{enumerate}
\end{theorem}

\pf

 $1\Rightarrow 2$:
 Write $G = \cap_{n=1}^{\infty} G_n$, where each $G_n$ is open.
 Fix $n$, and choose for each $x$ in $X$ an open interval $I^n_x$ which
 contains $x$, and whose closure is contained in $G_n$.
 Choose a countable subcover $\{I^n_{x^n_j}:j=1,2,3,...\}$ of $X$ from
 the cover $\{I^n_x:x \in X\}$. For each $n$ and $k$ define
 $I^n_k = \cup_{j\leq k}I^n_{x^n_j}$.

 Then ${\cal U}_n=\{I^n_k:k=1,2,3,...\}$ is a $\gamma$--cover of $X$
 such that for each $k$ the closure of $I^n_k$ is contained in $G_n$.

 Apply the fact that $X$ is a $\ufin(\ga,\ga)$-set to the
 sequence
 $$({\cal U}_n:n=1,2,3,...).$$
 For each $n$ choose a $k_n$ such that
 $(I^n_{k_n}:n=1,2,3,...)$  is a $\gamma$--cover of $X$.
 For each $n$ let
 $F_n$ be the intersection of the closures of the sets
 $I^m_{k_m},\ m\geq n$. For each $n$ we have the closed set $F_n$
 contained in $G$. But then the union of the $F_n$'s is an $F_{\sigma}$ 
 set which contains $X$ and is contained in $G$.

\medskip
 $2\Rightarrow 1$: Let $({\cal U}_n:n<\omega)$ be a sequence
 such that each ${\cal U}_n$ is a cover of $X$ by open subsets
 of the real line.   By assumption there
 exists closed sets $F_n$ such that
 $$X\subseteq \bigcup_{n<\omega} F_n\subseteq
 \bigcap _{n<\omega} (\cup {\cal U}_n).$$
 Since the real line is $\sigma$-compact we may assume that
 the $F_n$ are compact.   For each $n$ choose
 ${\cal V}_n\in [{\cal U}_n]^{<\omega}$ such that
 $(\cup_{m<n}F_m)\subseteq \cup {\cal V}_n$ for each $n$.
Either there exists
$n$ such that $\cup{\cal V}_n=X$ or
$\{\cup{\cal V}_n:n\in\omega\}$ is infinite and hence
a $\gamma$-cover of $X$.
\qed

\section{Ramseyan theorems and other properties} \label{ramsey}

Other classes of spaces motivated by diagonalization of open
covers are related to $Q$-point ultrafilters, $P$-point ultrafilters
and Ramsey-like partition relations. If ${\cal A}$ and ${\cal B}$
are classes of open covers, then a space has the property
\begin{enumerate}

\item ${\sf Q}({\cal A},{\cal B})$ iff for every open cover
${\cal U}\in {\cal A}$ and for every
partition of this cover into countably many disjoint nonempty finite
sets ${\cal F}_0,\ {\cal F}_1,\ {\cal F}_2,\dots$, there is a
subset ${\cal V}\subseteq{\cal U}$ which belongs to ${\cal B}$ such
that $|{\cal V}\cap{\cal F}_n|\leq 1$ for each $n$ and

\item ${\sf P}({\cal A},{\cal B})$ iff
for every sequence $\{{\cal U}_n: n\in\omega\}$ of open
covers of $X$ from ${\cal A}$ such that
${\cal U}_{n+1}\subseteq {\cal U}_n$, for each $n$,
there is an open cover ${\cal V}$ which belongs to
${\cal B}$ such that
${\cal V}\subseteq^*{\cal U}_n$ for each $n$.

\end{enumerate}

In Scheepers \cite{S} the partition relation
$\om\rightarrow(\om)^2_2$ was defined:
a space $X$ is said to satisfy $\om\rightarrow(\om)^2_2$ iff
for every $\omega$--cover ${\cal U}$ of $X$, if
  \[f:[{\cal U}]^2\rightarrow\{0,1\}\]
is any coloring, then there is an $i\in\{0,1\}$, an
$\omega$--cover ${\cal V}\subseteq {\cal U}$ such that
$f(\{A,B\})=i$
for all $A$ and $B$ from ${\cal V}$. It is customary to say
that ${\cal V}$ is homogeneous
for $f$.

Also in \cite{S} it was shown that for a set $X$ of real numbers,
the following statements are equivalent:
\begin{enumerate}
 \item{$X$ is both $\sone(\om,\om)$ and ${\sf Q}(\om,\om)$.}
 \item{$\om$, the collection of $\omega$--covers of $X$,
     satisfies the following partition
     relation: $\om\rightarrow(\om)^2_2$.}
\end{enumerate}
  The next theorem shows that indeed, the partition relation characterizes
the property of being a
   $\sone(\om,\om)$--set. This also implies that
 $$\sone(\om,\om)={\sf P}(\om,\om)+{\sf Q}(\om,\om).$$

\begin{theorem} $\sone(\om,\om)\subseteq {\sf Q}(\om,\om)$.
\end{theorem}
\pf Let $X$ be a $\sone(\om,\om)$--set and let ${\cal U}$ be an
$\omega$--cover of it.
   Let $({\cal P}_n:n<\omega)$ be a partition of this cover into
pairwise disjoint
   finite sets. Enumerate the cover bijectively as $(U_n:n<\omega)$
such that, letting
   for each $n$ the set $I_n$ be the $j$'s such that $U_j\in
{\cal P}_n$. We get a partition $(I_n:n<\omega)$ of $\omega$ into
disjoint intervals such
   that if $m$ is less than $n$, then each element of $I_m$ is
less than each
   element of $I_n$. For each $\ell$, let
$m_{\ell}=\sum_{j\leq\ell}|I_j|$.
  Now define an $\omega$--cover ${\cal V}$ of $X$ such that $V$ is
in ${\cal V}$ iff
\[V= U_{k_0}\cap\dots\cap U_{k_r}
\]
  where
\begin{enumerate}
  \item{$r=m_{\ell_0}$ and}
  \item{$\ell_0<\dots<\ell_r$ are such that for each $j$, $k_j$ is in
         $I_{\ell_j}$, and}
  \item{$V$ is nonempty.}
\end{enumerate}

   Next, choose a partition $({\cal V}_n:n<\omega)$ such that each
${\cal V}_n$ is an
   $\omega$--cover of $X$, and ${\cal V}$ is the union of these
sets. Then, discard from each
   ${\cal V}_n$ all sets of the form
\[U_{k_0}\cap\dots\cap U_{k_r}
\]
   where $k_0$ is an element of $I_0\cup\dots\cup I_n$; let
${\cal W}_n$ denote the resulting
family. Observe that each ${\cal W}_n$ is still an
$\omega$--cover.

   Since $X$ is an $\sone(\om,\om)$--set, we find for each $n$ a $W_n$ in
${\cal W}_n$ such that the
   set $\{W_n:n\in\omega\}$ is an $\omega$--cover of $X$. For
each $n$ we fix a representation
\[W_n = U_{k^n_0}\cap\dots\cap U_{k^n_{r(n)}}
\]
   where $k^n_0<\dots<k^n_{r(n)}$. On account of the way ${\cal
W}_n$ was obtained from ${\cal V}_n$,
   we see that $n<k^n_0$ and $n<r(n)$. Now choose recursively sets
$$U_{k(0)}, U_{k(1)},\dots,U_{k(n)},\dots$$
   so that $U_{k(0)} = U_{k^1_0}\supseteq W_1$. Suppose that
$U_{k(0)}, \dots, U_{k(n)}$ have been chosen such that for $i\leq
n$ we have
\begin{itemize}
\item{$k(i)\in\{k^i_0,\dots,k^i_{r(i)}\}$ and}
\item{$W_i\subseteq U_{k(i)}$, and}
\item{the $k(i)$'s belong to distinct $I_j$'s,}
\end{itemize}
To define $U_{k(n+1)}$ we consult $W_n=U_{k^{n+1}_0}\cap \dots\cap
U_{k^{n+1}_{r(n+1)}}$.
   Since we have so far selected only $n+1$ numbers and since
$r(n+1)$ is larger than $n+1$, and since the
   $k^{n+1}_j$ come from $r(n+1)$ disjoint intervals $I_j$, we can find
one of these intervals which is disjoint
   from $\{k(0),\dots,k(n)\}$, and select $k(n+1)$ to be the
$k^{n+1}_j$ from that interval. This then
   specifies $U_{k(n+1)}$.

   Because the sequence of $W_n$'s refines $\{U_{k(n)}:n<\omega\}$,
the latter is an $\omega$--cover of $X$,
   and by construction it contains no more than one element per
${\cal P}_n$.
\qed

   In Scheepers \cite{S} it was shown that if $X$ satisfies
$\sfin(\om,\om)$, then its family of $\omega$--covers,
$\om$, satisfies the partition relation
\[\om\rightarrow\lceil\om\rceil^2_2.
\]
Satisfying this partition relation means that for every
$\omega$--cover ${\cal U}$ of
$X$, if
\[f:[{\cal U}]^2\rightarrow\{0,1\}
\]
   is any coloring, then there is an $i\in\{0,1\}$, an
$\omega$--cover ${\cal V}\subseteq {\cal U}$ and a
   finite--to--one function $q:{\cal V}\rightarrow\omega$ such that
for all $A$ and $B$ from ${\cal V}$, if
   $q(A)\neq q(B)$, then $f(\{A,B\})=i$. It is customary to say
that ${\cal V}$ is eventually homogeneous
   for $f$.

  We now show that these two properties are equivalent.

\begin{theorem} For any space $X$,
$\om\rightarrow\lceil\om\rceil^2_2$ is equivalent
to $\sfin(\om,\om)$.
\end{theorem}

\noindent {\pf}
$\sfin(\om,\om)$ implies 
$\om\rightarrow\lceil\om\rceil^2_2$
is proved in \cite{S}.
To prove the other direction suppose that
${\cal U}_{n}=\{ U^n_m:m\in\omega\}$ is an
$\omega$--cover for each $n\in \omega$. Let
$$
{\cal U}=\{ U^0_k\cap U^k_l:k,l\in\omega\}.
$$
${\cal U}$ is an $\omega$--cover, since given a finite
$F\subseteq X$ we can first pick $k$ with $F\subseteq U^0_k$ and
then pick $l$ with $F\subseteq U^k_l$.  For each
element of ${\cal U}$ we pick a pair as above and
define $f:[{\cal U}]^{2}\rightarrow 2$ by
$$f(\{U^0_{k_{0}}\cap U^{k_{0}}_{l_{0}},
U^0_{k_{1}}\cap U^{k_{1}}_{l_{1}}\})=
\left\{
\begin{array}{ll}
0 & \mbox{if $k_0=k_1$}\\
1 & \mbox{if $k_0\not=k_1$}
\end{array}
\right.$$

By applying $\om\rightarrow\lceil\om\rceil^2_2$
there exists a sequence
$(k_{i},l_{i})$ and a finite--to--one function
$q:\omega\rightarrow\omega$ such that
$${\cal V}=\{U^0_{k_{i}}\cap U^{k_{i}}_{l_{i}}: i \in\omega\}$$
is
an $\omega$--cover of $X$ and either
\begin{itemize}
\item[(a)]{$q(i)\not=q(j)$ implies $k_i=k_j$ or }
\item[(b)]{$q(i)\not=q(j)$ implies $k_i\not=k_j$. }
\end{itemize}

In case (a) since $q$ is  finite--to--one,
we get that $k_i=k_j$ for every $i,j\in \omega$.  This would
mean that every element
of ${\cal V}$ refines $U^0_{k_{0}}$,  but this
contradicts the fact that ${\cal V}$ is an $\omega$--cover.
Thus this case cannot occur.

In case (b) let
$${\cal W}=\{U^{k_i}_{l_i}: i<\omega\}.$$
Since ${\cal V}$ refines ${\cal W}$ and
$X\notin {\cal W}$,
${\cal W}$ is an $\omega$--cover of $X$.
Define
$${\cal W}_n=\{U^{k_i}_{l_i}: k_i=n\}\subseteq {\cal U}_n.$$
To finish
the proof it is enough to see that each ${\cal W}_n$ is
finite.   If not, there would be an infinite
$A\subseteq\omega$ such that $k_i=n$ for each $i\in A$.
Since $q$ is finite--to--one,
there would be $i\not=j\in A$ with $q(i)\not=q(j)$.  But
$k_i=k_j=n$ contradicts the assumption of case (b).
\qed

\bigskip

\begin{center}
  Addresses
\end{center}

\medskip

\begin{flushleft}
Winfried Just                             \\
Ohio University                           \\
Department of Mathematics                 \\
Athens, OH 45701-2979 USA                 \\
e-mail: just@ace.cs.ohiou.edu             \\
\end{flushleft}

\begin{flushleft}
 Arnold W. Miller                                     \\
 University of Wisconsin-Madison                      \\
 Department of Mathematics  Van Vleck Hall            \\
 480 Lincoln Drive                                    \\
 Madison, Wisconsin 53706-1388, USA                   \\
 e-mail: miller@math.wisc.edu                         \\
 home page: http://math.wisc.edu/$^\sim$miller/      \\
\end{flushleft}

\begin{flushleft}
Marion Scheepers                    \\
Department of Mathematics           \\
Boise State University              \\
Boise, Idaho 83725 USA              \\
e-mail: marion@math.idbsu.edu       \\
\end{flushleft}

\begin{flushleft}
Paul J. Szeptycki                         \\
Ohio University                           \\
Department of Mathematics                 \\
Athens, OH 45701-2979 USA                 \\
e-mail:  szeptyck@oucsace.cs.ohiou.edu    \\
\end{flushleft}

\bigskip

\begin{center}
   August 1995
\end{center}

\end{document}